\documentclass[review]{elsarticle}
	\usepackage{lineno,hyperref}
	\journal{Journal of \LaTeX\ Templates}
	\usepackage[russian,english]{babel}
	\usepackage{amsthm}
	\theoremstyle{plain}

	\newtheorem{theorem}{Theorem}
	\newtheorem{definition}{Definition}
	\newtheorem{lemma}{Lemma}
	\newtheorem{corollary}{Corollary}
	\newtheorem*{proof*}{Proof}
	\usepackage{lineno,hyperref}
	\usepackage{ucs} 
	\usepackage[russian]{babel}  % Включаем пакет для поддержки русского языка 
	\usepackage{tikz} 
	
\begin{document}

		\begin{frontmatter}
			
			\title{Poincare--Riemann--Hilbert boundary-value problem for The  Millennium Prize Problems .}
			%\tnotetext[mytitlenote]{Fully documented templates are available in the elsarticle package on %\href{http://www.ctan.org/tex-archive/macros/latex/contrib/elsarticle}{CTAN}.}
			
			\author{Durmagambetov A.A\fnref{myfootnote}}
			\address{}
			\fntext[myfootnote]{}
			\ead[url]{}
			
			\ead{aset.durmagambet@gmail.com}
			
			\address[mymainaddress]{010000, Kazakhstan }
			
			\begin{abstract}
				Using the example of a complicated problem such as the Cauchy problem for the Navier--Stokes equation, we show how the Poincar\'e--Riemann--Hilbert boundary-value problem enables us to construct effective estimates of solutions for this case. The apparatus of the three-dimensional inverse problem of quantum scattering theory is developed for this. It is shown that the unitary scattering operator can be studied as a solution of the Poincar\'e--Riemann--Hilbert boundary-value problem.  The same scheme of reduction of Riemann integral equations for the zeta function to the Poincar\'e--Riemann--Hilbert boundary-value problem allows us to construct effective estimates that describe the behaviour of the zeros of the zeta function very well.
			\end{abstract}
			
			
		\end{frontmatter}

%	\maketitle
	
	\section {Introduction}
	Using the example of a complicated problem such as the Cauchy problem for the Navier--Stokes equation, we show how the Poincar\'e--Riemann--Hilbert boundary-value problem enables us to construct effective estimates of solutions for this case. The apparatus of the three-dimensional inverse problem of quantum scattering theory is developed for this. It is shown that the unitary scattering operator can be studied as a solution of the Poincar\'e--Riemann--Hilbert boundary-value problem. This allows us to go on to study the potential in the Schr\"odinger equation, which we consider as a velocity component in the Navier--Stokes equation. The same scheme of reduction of Riemann integral equations for the zeta function to the Poincar\'e--Riemann--Hilbert boundary-value problem allows us to construct effective estimates that describe the behaviour of the zeros of the zeta function very well.

	\section{Results for the one-dimensional case}
	%\subsection{General}
	Let us consider a one-dimensional function $ {f}  $ and its Fourier transformation $  \tilde{f} $. Using the notions of module and phase, we write the Fourier transformation in the following form: $  \tilde{f}=|\tilde{f}|\exp(i\Psi) $ , where $ \Psi $ is the phase.  The Plancherel equality states that $ ||f||_{L_{2}}={\rm const}||\tilde{f} ||_ {L_{2}} $. Here we can see that the phase does not contribute to determination of the $X$ norm. To estimate the maximum we make a simple estimate as $  {\rm max}|f|^2 \le 2||f||_{L_{2}} ||\nabla f||_{L_{2}}   $. Now we have an estimate of the function maximum in which the phase is not involved. Let us consider the behaviour of a progressing wave travelling with a constant velocity of $  v=a$ described by the function $ {F(x,t)=f(x+at)}  $.  Its Fourier transformation with respect to the variable $x$ is $  \tilde{F}=\tilde{f}{\exp}(iatk) $.  Again, in this case, we can see that when we study a module of the Fourier transformation, we will not obtain major physical information about the wave, such as its velocity and location of the wave crest because $ |\tilde{F}|=|\tilde{f}| $ . These two examples show the weaknesses of studying the Fourier transformation. Many researchers focus on the study of functions using the embedding theorem, in which the main object of the study is the module of the function. However, as we have seen in the given examples, the phase is a principal physical characteristic of any process, and as we can see in mathematical studies that use the embedding theorem with energy estimates, the phase disappears. Along with the phase, all reasonable information about the physical process disappears, as demonstrated by Tao [1] and other research studies. In fact, Tao built progressing waves that  are  not followed by energy estimates . Let us proceed with a more essential analysis of the influence of the phase on the behaviour of functions.
	\begin  {theorem}
	There are functions of $ W_{2}^{1}(R)$ with a constant rate of the norm for a gradient catastrophe for which a phase change of its Fourier transformation is sufficient.
	\end  {theorem}
	Proof: To prove this, we consider a sequence of testing functions $  \tilde{f_{n}}=\Delta/(1+k^2),\;\Delta=(i-k)^n/(i+k)^n $. It is obvious that $|\tilde{f_{n}}|=1/(1+k^2)$ and
	$  {\rm max}|f_{n}|^2 \le 2||f_{n}||_{L_{2}} ||\nabla f_{n}||_{L_{2}} \le {\rm const}   $. Calculating the Fourier transformation of these testing functions, we obtain 
	\begin{equation}
	f_{n}(x) =	x(-1)^{(n-1)}2\pi \exp(-x)L^1_{(n-1)}(2x)  {\rm if}\,\, x> 0 ,\,\,
	f_{n}(x) =		0  \,\,\,if \,\,\, x\le 0, 
	\end{equation}
	where $ L^1_{(n-1)}(2x)$ is a Laguerre polynomial.
	Now we see that the functions are equibounded and derivatives of these functions will grow with the growth of $ {n}.  $
	Thus, we have built an example of a sequence of the bounded functions of $ W_{2}^{1}(R)$ which have a constant norm  $ W_{2}^{1}(R)$, and this sequence converges to a discontinuous function.\\
	The results show the flaws of the embedding theorems when analyzing the behavior of functions. Therefore, this work is devoted to overcoming them and the basis for solving the formulated problem is the analytical properties of the Fourier transforms of functions on compact sets. Analytical properties and estimates of the Fourier transform of functions are studied using the Poincaré – Riemann – Hilbert boundary value problem
	
	%------------------------------------------------
	\section{Results for the three-dimensional case }
	Consider Schr\"odinger's equation:
	\begin{equation}
	\label {eq:1}
	-\Delta_{x}\Psi  +q\Psi  =k^{2}\Psi, ~k\in C.  
	\end{equation}
	Let $ \Psi_{+ } (k,\theta,x)  $ be a solution of (\ref{eq:1}) with the following asymptotic behaviour:
	\begin{equation}\label{eq:2}
	\Psi_{+ } (k,\theta,x)=\Psi _{0}(k,    
	\theta, x)+\frac{e^{ ik|x|}}{|x|}A(k,\theta^{'},\theta)+0\left( \frac{1}{|x|}\right), \; |x|\rightarrow \infty,
	\end{equation}
	where   $ A(k,\theta^{'},\theta)  $ is the scattering amplitude and  $ \theta^{'}= \frac{x}{|x|},\;\theta\in S^{2} $ for  $ k\in \bar{C}^{+}=\{{\rm Im} k\ge 0\} $ $ \Psi _{0}(k,\theta, x) =e^{ik(\theta, x)}  $:
	$$\label {eq:3}
	A(k,\theta^{'},\theta)=-\frac{1}{4\pi }\int_{R^{3}}q(x)\Psi _{+ } (k,\theta,x)e^{-ik\theta^{'} x}dx.
	$$
	Solutions to (\ref{eq:1}) and (\ref{eq:2}) are obtained by solving the integral equation
	$$\label{eq:4}
	\Psi_{+ } (k,\theta,x)=\Psi _{0}(k,    
	\theta, x)+\int_{R^3}q(y)\frac{e^{+ ik|x-y|}}{|x-y|}\Psi{+ } (k,\theta,y)dy=G(q\Psi_{+ }), 
	$$
	which is called the Lippman--Schwinger equation.
	
	Let us introduce
	$$ 
	\theta,\theta^{'}\in S^{2}, 
	Df=k\int_{S^{2}}A(k,\theta^{'},\theta)f(k,\theta^{'}) d\theta^{'}.
	$$

	% % % % % % % % % % % % % % % % % % % % % % % % % % % % % % % % % % % % % % % % % % % % % % % % % % % % % % % % % % % % % % % % % %
	Let us also define the solution $ \Psi _{- } (k,\theta,x)$ for $ k\in \bar{C}^{-}=\{{\rm Im} k\le 0\} $ as  \[\Psi _{- } (k,\theta,x)=\Psi _{+ } (-k,-\theta,x) .\]
	As is well known [8], %\cite{a1}
	\begin{equation}
	\label{eq:eq5}
	\Psi _{+ } (k,\theta,x)-\Psi _{- } (k,\theta,x)=
	-\frac{k}{4\pi}\int_{S^{2}}A(k,\theta^{'},\theta)\Psi_{-} (k,\theta^{'},x)d\theta^{'},~  k\in R.
	\end{equation}
	This equation is the key to solving the inverse scattering problem and was first used by Newton~[8,9] and Somersalo et al.~[10].
	
	\begin{definition}
		The set of measurable functions $\mathbf{R}$ with the norm defined by $$
		||q||_{\mathbf{R}}=\int_{R^{6}}\frac{q(x)q(y)}{|x-y|^{2}}dxdy<\infty $$
		is recognised as being of Rollnik class.
		
	\end{definition}\label{df:d1}

	Equation (\ref{eq:eq5}) is equivalent to the following:
	$$
	\Psi _{+}=S\Psi _{-},  \label{eq:eq6}
	$$
	where  $ S $ is a scattering operator with the kernel 
	$$ ~S(k,\textit{\l})=\int_{R^{3}}\Psi _{+}(k,x)\Psi _{-}^{\ast }(\textit{\l},x)dx. $$
	
	The following theorem was stated in  [9]:
	
	\begin{theorem}\label{thm:t1}\textbf{(Energy and momentum conservation laws)}
	Let $q\in \mathbf{R}$. Then, $ SS^{\ast }=I$ and $S^{\ast }S=I,$ where $ I $ is a unitary operator.
	\end{theorem}
	\begin{corollary}  $ SS^{\ast }=I$ and $S^{\ast }S=I$  yield
	$$
	A(k,\theta^{'},\theta)-A(k,\theta,\theta^{'})^{\ast }=\frac{ik}{2\pi}\int_{S^{2}}A(k,\theta,\theta^{''})A (k,\theta^{'}, \theta^{''} )^{\ast }d\theta^{''}.
	$$
	
	\end{corollary}
	
	\begin{theorem}
	\label{Theorem 1.4}
	\textbf{(Birmann--Schwinger estimation)}
	Let $q\in \mathbf{R}$. Then, the number of discrete eigenvalues can be estimated as 
	$$
	N(q)\leq \frac{1}{(4\pi )^{2}}\int_{R^{3}}\int_{R^{3}}\frac{%
		q(x)q(y)}{|x-y|^{2}}dxdy.
	$$
	\end{theorem}
	%%%
	\begin{lemma}\label{lm:l2} %\Delta=\prod\limits_{j=1}^{N} (k+iE_{j})/(k-iE_{j}
	Let $
	\left(|q|_{L_{1}(R^3)} +4\pi|q|_{L_{2}(R^3)}   \right)<\alpha<1/2   $. Then, 
	$$
	\left\| \Psi_{+}\right\|_{L_{\infty}} \le \frac{\left(|q|_{L_{1}(R^3)} +4\pi|q|_{L_{2}(R^3)}   \right)}{ 1-\left(|q|_{L_{1}(R^3)} +4\pi|q|_{L_{2}(R^3)}   \right)    } < \frac{\alpha}{1-\alpha}, 
	$$
	$$
	\left\|   \frac{ \partial(\Psi_{+}-\Psi_0)}{\partial k}\right\|_{L_{\infty}}\le  \frac{  	|q|_{L_{1}(R^3)} +4\pi|q|_{L_{2}(R^3)}  }{1-\left(|q|_{L_{1}(R^3)} +4\pi|q|_{L_{2}(R^3)}   \right)    }  < \frac{\alpha}{1-\alpha}     
	.$$
	\end{lemma}
	\begin{proof}
	By the Lippman--Schwinger equation, we have
	$$
	\left| \Psi_{+}-\Psi_0\right| \le  \left| Gq\Psi_{+}\right|, \,\,\,\,
	$$
	
	$$
	\left| \Psi_{+}-\Psi_0\right|_{L_{\infty}} \le \left| \Psi_{+}-\Psi_0\right|_{L_{\infty}} \left|   Gq\right|+\left|Gq\right| 	,$$
	and, finally,	
	$$
	\left| \Psi_{+}-\Psi_0\right| \le \frac{\left(|q|_{L_{1}(R^3)} +4\pi|q|_{L_{2}(R^3)}   \right)}{ 1-\left(|q|_{L_{1}(R^3)} +4\pi|q|_{L_{2}(R^3)}   \right)    }.  
	$$	
	
	By the Lippman--Schwinger equation, we also have
	$$
	\left| \frac{ \partial\left( \Psi_{+}-\Psi_0\right) }{\partial k} \right| \le  \left|\frac{ \partial Gq}{\partial k}   \Psi_{+}\right|+  \left| Gq \frac{ \partial\left( \Psi_{+}-\Psi_0\right) }{\partial k}   \right| +\left| Gq\right|,
	$$
	$$\label{eq:eq2}	\left| \frac{ \partial (\Psi_{+}-\Psi_0)}{\partial k} \right|  \le \left(|q|_{L_{1}(R^3)} +4\pi|q|_{L_{2}(R^3)}   \right),  
	$$
	
	$$
	\left\|   \frac{ \partial(\Psi_{+}-\Psi_0)}{\partial k}\right\|_{L_{\infty}}\le  \frac{  	|q|_{L_{1}(R^3)} +4\pi|q|_{L_{2}(R^3)}  }{1-\left(|q|_{L_{1}(R^3)} +4\pi|q|_{L_{2}(R^3)}   \right)    },    
	$$
	which completes the proof.
	\end{proof}
	Let us introduce the following notation:
	$$ 
	Q(k,\theta,\theta^{'})=\int_{R^3}q(x)e^{ik(\theta-\theta^{'})x}dx,\;K(s)=s,\;X(x)=x,\,\,\, 	$$
	$$\label{eq:eq2} T_+Q=\int_{-\infty}^{+\infty}\frac{Q(s,\theta,\theta^{'})}{ s-t-i0}ds,\,\,\,T_-Q=\int_{-\infty}^{+\infty}\frac{Q(s,\theta,\theta^{'})}{ s-t+i0}ds .
	$$

	\begin{lemma}\label{lm:l2} %\Delta=\prod\limits_{j=1}^{N} (k+iE_{j})/(k-iE_{j}
	Let $q\in \mathbf{R}\cap  L_{1}(R^3),\; \left\| q \right\|_{L_1}+4\pi|q|_{L_{2}(R^3)} <\alpha<1/2   $. Then, 
	$$
	\left\|  	A _{+}\right\|_{L_{\infty}}  < \alpha+  \frac{\alpha}{1-\alpha} ,  
	$$
	
	$$
	\left\|  \frac{ \partial A _{+}}{\partial k}	\right\|_{L_{\infty}}  < \alpha+  \frac{\alpha}{1-\alpha}.
	$$
	
	\end{lemma}
	\begin{proof}
	Multiplying the Lippman--Schwinger equation by $q(x)\Psi_0(k,\theta,x)$ and then integrating, we have
	$$
	A(k,\theta,\theta^{'})=	Q(k,\theta,\theta^{'})+		\int_{R^{3}}q(x)\Psi_0(k,\theta,x)G	q\Psi_{+}dx.
	$$
	We can estimate this latest equation as
	$$
	\left| A\right| \le \alpha+ \alpha\frac{\left(|q|_{L_{1}(R^3)} +4\pi|q|_{L_{2}(R^3)}   \right)}{ 1-\left(|q|_{L_{1}(R^3)} +4\pi|q|_{L_{2}(R^3)}   \right)    }.
	$$
	Following a similar procedure for $	\left\|  \frac{ \partial A _{+}}{\partial k}	\right\|$ completes the	proof.	
	\end{proof}
	
	We define the operators ~$T_{\pm }$, $T$ for ~$f\in W_{2}^{1}(R)$ as follows:
	$$ 
	T_{+}f=\frac{1}{2\pi i}\lim\limits_{{\rm Im}z\rightarrow 0}\int\limits_{-\infty}^{\infty }\frac{f({s})}{s-z}ds,~{\rm Im}~z>0,~
	T_{-}f=\frac{1}{2\pi i}\lim\limits_{{\rm Im}z\rightarrow 0}\int\limits_{-\infty
	}^{\infty }\frac{f({s})}{s-z}ds,~{\rm Im}~z<0,
	$$
	$$
	Tf=\frac{1}{2}(T_{+}+T_{-})f.
	$$
	Consider the Riemann problem of finding a function $\Phi $ that is analytic in the complex plane with a cut along the real axis. Values of
	$ \Phi $ on the two sides of the cut are denoted as $\Phi_{+}$ and $\Phi_{-} $. The following presents the results of [12]:
	\begin{lemma}\label{lm:l1}
	$$
	TT= \frac14 I,~TT_+ = \frac12 T_+,~TT_- = - \frac12 T_-, \ T_+ = T+\frac12 I,~T_- = T-\frac12 I,~T_-T_- = -  T_-
	.$$ 
	\end{lemma}
	
	Denote
	$$
	\Phi_{+}(k,\theta,x)=\Psi_{+}(k,\theta,x)-\Psi_{0}(k,\theta,x), \,\,\,\Phi_{-}(k,\theta,x)=\Psi_{-}(k,-\theta,x)-\Psi_{0}(k,\theta,x),\,\, 	$$
	$$\label{eq:eq2} g(k,\theta,x)=\Phi_{+}(k,\theta,x)-\Phi_{-}(k,\theta,x)/
	$$
	
	%\Phi_{-}(k,\theta,x)=\Psi_{-}(k,-\theta,x)-\ph0(k,\theta,x)
	
	\begin{lemma}\label{lm:l2} %\Delta=\prod\limits_{j=1}^{N} (k+iE_{j})/(k-iE_{j}
	Let $q\in \mathbf{R},\,\, N(q)<1, ~g_{+}=g(k,\theta,x)$, and $g_{-}=g(k,-\theta,x). $ Then, 
	$$
	\Phi _{+}(k,\theta,x)=T_{+ } g_{+}  +e^{ik\theta x},\ \Phi _{-}(k,\theta,x)=T_{- } g_{+}  +e^{ik\theta x}.
	$$ 
	\end{lemma}	
	\begin{proof}  
	The proof of the above follows from the classic results for the Riemann problem.
	\end{proof}
	\begin{lemma}\label{lm:l2} %\Delta=\prod\limits_{j=1}^{N} (k+iE_{j})/(k-iE_{j}
	Let $q\in \mathbf{R},\,\, N(q)<1, ~g_{+}=g(k,\theta,x),$ and $g_{-}=g(k,-\theta,x) , ) $. Then, 
	$$
	\Psi _{+}(k,\theta,x)=(T_{+ } g_{+}  +e^{ik\theta x}),\ \Psi _{-}(k,\theta,x)=(T_{- } g_{-}  +e^{-ik\theta x}).
	$$ 
	\end{lemma}
	\begin{proof}
	The proof of the above follows from the definitions of $ g,$ $\Phi _{\pm }$, and $\Psi _{\pm } $ .
	\end{proof}
	%%%%%%%%%%%%%%%%%%%%%%%%%%%%%%%%%%%%%%%
	
	\begin{lemma}\label{lm:l6}
	Let $$ \sup\limits_{ k} \left|    \int\limits_{-\infty}^{\infty} \frac{ pA(p,\theta^{'},\theta) }{ 4\pi(p-k+i0 ) }dp  \right|<\alpha,\,\,\int_{S_2}\alpha d\theta<1/2	.$$ %$$	 \sup\limits_{ k} \left|    \int\limits_{-\infty}^{\infty} \frac{ pA(p,\theta^{'},\theta) \Psi_{0}}{ 4\pi(p-k+i0 ) }dp  \right|<\alpha<1/2 
	%	$$ 
	Then,
	$$
	\prod\limits_{ 0\le j<n }  \int_{S_2}\left|   \int_{-\infty}^{\infty} \frac{ {k_{j}}A(k_{j},\theta^{'}_{k_{j}},\theta_{k_{j}}) }{ 4\pi( k_{j+1}-k_{j}+i0 ) }d{k_{j}} \right |d\theta_{k_{j}} \le 2^{-n}.
	\label{lm:psi}
	$$
	\end{lemma}
	\begin{proof}
	
	Denote 
	$$
	\alpha_j{}=\left| Vp \int_{-\infty}^{\infty} \frac{ {k_{j}}A(k_{j},\theta^{'}_{k_{j}},\theta_{k_{j}}) }{ 4\pi( k_{j+1}-k_{j}+i0 ) }d{k_{j}}\right| , \,\,\, 
	%\beta_{j} = \frac{1}{2}  {k_{j}}A(k_{j},\theta^{'}_{k_{j}},\theta_{k_{j}}).
	$$
	Therefore,
	$$
	\prod\limits_{ 0\le j<n } \int_{S_2} \left| \int_{-\infty}^{\infty} \frac{ {k_{j}}A(k_{j},\theta^{'}_{k_{j}},\theta_{k_{j}}) }{ 4\pi( k_{j+1}-k_{j}+i0 ) }d{k_{j}}  \right |d\theta_{k_{j}}\le   \prod\limits_{ 0\le j<n } \int_{S_2} \alpha_j{} d\theta_{k_{j}}  <
	2^{-n}.	$$
	$$\label{eq:eq2}
	%   as\,\,\, m\neq l ,\alpha_{j_{m}}\neq \alpha_{j_{l}} \,\,\,and \,\,\beta_{j_{m}}\neq \beta_{j_{l}}
	$$
	This completes the proof.
	\end{proof}
	%%%%%%%%%%%%%%%%%%%%%%%%%%%%%%%%%%%%%%%%%%%%%%%%%%%%%%%%%%%%%%%%%%%

	%%%%%%%%%%%%%%%%%%%%%%%%%%%%%%%%%%%%%%%

	\begin{lemma}\label{lm:TA}
	
	Let 	$$
	\sup\limits_{ k } \int_{S^{2}} \left| T_{-}QK\right|d\theta \le \alpha<\frac{1}{2C}<1,\,\,\,\, \sup\limits_{ k } \int_{S^{2}}\left| T_{-}\tilde{q}K\right|d\theta\le \alpha<\frac{1}{2C}<1	,\,\,\,\, 	$$
	$$\label{eq:eq2}\sup\limits_{ k } \int_{S^{2}}\left| T_{-}Q\tilde{q}K^2\right|d\theta\le \alpha<\frac{1}{2C}<1.		$$
	Then,
	$$
	\sup\limits_{ k}\int_{S^{2}}\left| 	T_{-}AK\right|d\theta \le\frac{ C\int_{S^{2}}\left|T_{-}QK\right|d\theta}{1-\sup\limits_{ k}\int_{S^{2}}\left| 	T_{-}A\tilde{q}K^2\right|d\theta},\,\,\, 		$$
	$$\label{eq:eq2}	
	\sup\limits_{ k}\left| \int_{S^{2}}	T_{-}A\tilde{q}K^2d\theta\right| \le
	\frac{ C\left|T_{-}\int_{S^{2}}Q\tilde{q}K^2d\theta\right|}{1- \left| T_{-}\int_{S^{2}}\tilde{q}Kd\theta\right|	}.
	$$

	\end{lemma}
	\begin{proof}
	By the definition of the amplitude and Lemma 4, we have 
	$$
	A(k,\theta^{'},\theta)=-\frac{1}{4\pi }\int_{R^{3}}q(x)\Psi _{+ } (k,\theta,x)e^{-ik\theta^{'} x}dx 	$$
	$$\label{eq:eq2}=-\frac{1}{4\pi }\int_{R^{3}}q(x) \left[ e^{ik\theta^{'} x} +T_+g(k,\theta,\theta^{'}) \right] e^{-ik\theta^{'} x}dx   .
	$$%
	We can rewrite this as
	\begin{equation}\label{eq:5}
	A(k,\theta^{'},\theta)=-\frac{1}{4\pi }\int_{R^{3}}q(x) \left[ e^{ik\theta x} +\sum_{n\ge 0}(-T_-D)^n\Psi_0\right]   e^{-ik\theta^{'} x} dx.
	\end{equation}
	Lemma 6 yields 
	$$
	\sup\limits_{ k}\int_{S^{2}}\left| 	T_{-}AK\right|d\theta\le	\sup\limits_{ k}	\int_{S^{2}}\left| \frac{1}{4\pi }T_{-}QK\right|d\theta  +\frac{\left( \sup\limits_{ k}\int_{S^{2}}\left| 	T_{-}KA\right|d\theta\right) ^2\int_{S^{2}}\left|T_{-} A\tilde{q}K^{2}\right|d\theta }{\left( 1-\sup\limits_{ k}\int_{S^{2}}\left| 	T_{-}KA\right|d\theta\right) ^{2}}. 
	$$
	
	Owing to the smallness of the terms on the right-hand side, the following estimate follows:
	$$
	\sup\limits_{ k}\int_{S^{2}}\left| 	T_{-}AK\right|d\theta\le	2\sup\limits_{ k}	\int_{S^{2}}\left| \frac{1}{4\pi }T_{-}QK\right|d\theta . 
	$$
	Similarly,
	$$\label{eq:eq2}		\sup\limits_{ k}\int_{S^{2}}\left| 	T_{-}A\tilde{q}K^2\right|d\theta \le C\int_{S^{2}}\left|T_{-}Q\tilde{q}K^2\right|d\theta+  \int_{S^{2}}\left| T_{-}A\tilde{q}K^2\right|d\theta	\int_{S^{2}} \left| T_{-}\tilde{q}K\right|d\theta,	   
	$$	
	$$\label{eq:eq2}	
	\sup\limits_{ k}\int_{S^{2}}\left| 	T_{-}A\tilde{q}K^2\right|d\theta \le
	\frac{ C\int_{S^{2}}\left|T_{-}Q\tilde{q}K^2\right|d\theta}{1- \int_{S^{2}}\left| T_{-}\tilde{q}K\right|d\theta	},
	$$
	$$
	\sup\limits_{ k}\int_{S^{2}}\left| 	T_{-}A\tilde{q}K^2\right|d\theta\le	2\sup\limits_{ k}\int_{S^{2}}	\left| \frac{1}{4\pi }T_{-}Q\tilde{q}K^2\right|d\theta  .
	$$

	This completes the proof.	
	\end{proof}
	%%%%%%%%%%%%%%%%%%%%%%%%%%%%%%%%%%%%%%%%%%%%%%%%%%%%%%%%%%%%%%%%%%%%%
	
	%%%%%%%%%%%%%%%%%%%%%%%%%%%%%%%%%%%%%%%%

	%%%%%%%%%%%%%%%%%%%%%%%%%%%%%%%%%%%%%%%%%%%%%%%%%%%%%%%%%%%%%%%%%%%%%%%%%%%%%%%%%%%%%%%%%%%%%%%%%
	To simplify the writing of the following calculations, we introduce the set defined by
	$$M_{\epsilon}(k)=\left( s|  \epsilon<|s|+|k-s|<\frac{1}{\epsilon}\right) .$$
	The Heaviside function is given by 
	$$
	{\Theta}(x) =
	\left\{
	1 ,\,\,\mbox{if } x> 0 ,\,\,\,
	 \,\,\,\,
	-1 \,\,  \mbox{if } x< 0 \,\,\,
	\right\}.
	$$

	\begin{lemma} \label{lm:l8}
	Let $q,\nabla q\in \cap L_{2}(R^{3})$, $|A|>0$. Then,
	$$ 
	\pi i\int_{R^3} \Theta(A) e^{ik|x|A}q(x)dx  =    \lim\limits_{\epsilon\rightarrow 0}\int_{s\in M_{\epsilon}(k)}\int_{R^3} \frac{e^{is|x|A}}{k-s} q(x)dxds,   
	$$
	$$ 
\pi i	\int_{R^3} \Theta(A) ke^{ik|x|A}q(x)dx  =    \lim\limits_{\epsilon\rightarrow 0}\int_{s\in M_{\epsilon}(k)}\int_{R^3} s\frac{e^{is|x|A}}{k-s} q(x)dxds  . 
	$$
	
	\end{lemma} 
	\begin{proof}
	The lemma can be proved by the conditions of lemma and the lemma of Jordan.
	\end{proof}

	%%%%%%%%%%%%%%%%%%%%%%%%%%%

	%%%%%%%%%%%%%%%%%%%%%%

	%%%%%%%%%%%%%%%%%%%%%%%%%%%%%%%%%%%%%%%%%%%%%%%%%%%%%%
	\begin{lemma}\label{lm:TkQ}
	
	Let $$ l=2	,\,\,\,I_0= \Psi_0(x,k)|_{ r=r_{0}} .
	$$
	Then
	$$		\left| 	\int_{-\infty}^{+\infty}\int_{S^{2}}\int_{S^{2}}\tilde{q}(k(\theta-\theta'))I_0k^2dkd\theta d\theta'\right| \le \sup\limits_{x\in R^3} \left| q(x)\right| +C_0(\frac{1}{r_{0}} +r_0)\left\| q\right\|_{L_2(R^3)}, 
	$$
	
	$$	 	\sup \limits_{\theta\in S^2}\left|\int_{-\infty}^{+\infty} \int_{S^{2}}\int_{S^{2}} QTKQI_0k^2d\theta''d\theta'dk \right|  \le  C_0(\frac{1}{r_{0}} +r_0)\left\| q\right\|^2_{L_2(R^3)} .
	$$	
	\end{lemma}
	\begin{proof}
	By the definition of the Fourier transform, we have
	$$\label{eq:eq2}	
	\int_{-\infty}^{+\infty}\int_{S^{2}}\int_{S^{2}}\tilde{q}(k(\theta-\theta'))I_0k^2dkd\theta d\theta'=
	\int_{-\infty}^{+\infty}\int_{S^{2}}\int_{S^{2}} \int_{0}^{+\infty}q(x)e^{ikx(\theta-\theta')}e^{ix_0k}k^2dkd\theta d\theta' drd\gamma,
	$$	
	where $x=r\gamma$
	The lemma of Jordan completes the proof for the first inequality.
	The second inequality is proved like the first:
	$$	 \int_{-\infty}^{+\infty} \int_{S^{2}}\int_{S^{2}} QTKQI_0k^2d\theta''d\theta'dk 
	$$
	$$
	=\int_{-\infty}^{+\infty} \int_{-\infty}^{+\infty} \int_{S^{2}}\int_{S^{2}}\int_{S^{2}} \frac{\left( \tilde{q}(s\cos(\theta')- s\cos(\theta''))\tilde{q}(k\cos(\theta)-s\cos(\theta'')\right)s }{k-s}   I_0k^2d\theta'd\theta''d\theta dk ds.
	$$				
	Lemma 8 yields			
	$$	\int_{-\infty}^{+\infty}\int_{S^{2}} \int_{S^{2}}\int_{S^{2}} \left( \tilde{q}(k\cos(\theta')- k\cos(\theta))\tilde{q}(k\cos(\theta)-k\cos(\theta'')\right) I_0k^3\Theta (\cos(\theta''))d\theta'd\theta''d\theta dk-
	$$
	$$	\int_{-\infty}^{+\infty} \int_{S^{2}}\int_{S^{2}} \int_{S^{2}}\left( \tilde{q}(k\cos(\theta')- k\cos(\theta))\tilde{q}(k\cos(\theta)-k\cos(\theta'')\right) I_0k^3\Theta (-\cos(\theta''))d\theta'd\theta''d\theta dk.
	$$
	Integrating $\theta$, $\theta'$, $\theta''$, and $k$, we obtain the proof of the second inequality of the lemma.
	
\end{proof}

%%%%%%%%%%%%%%%%%%%%%%%%%%%%%%%%%%%%%%

\begin{lemma}\label{lm:TA}
	
	Let 	$$
	\sup\limits_{ k }  \left| T_{-}QK\right| \le \alpha<\frac{1}{2C}<1,\,\,\,\, \sup\limits_{ k } \left| T_{-}\tilde{q}K\right|\le \alpha<\frac{1}{2C}<1	,\,\,\,\, 	$$
	$$\label{eq:eq2}\sup\limits_{ k } \left| T_{-}Q\tilde{q}K^2\right|\le \alpha<\frac{1}{2C}<1	,\,\,\,l=0,1,2	.	$$
	Then,
	
	$$\label{eq:A^0}	
	\left| 	\int_{-\infty}^{+\infty} \int_{S^{2}}\int_{S^{2}} A(k,\theta',\theta)k^ldkd\theta'd\theta \right| \le \left| \int_{-\infty}^{+\infty}\int_{S^{2}}\int_{S^{2}}\tilde{q}(k(\theta-\theta'))k^ldkd\theta' d\theta\right| 
	$$$$
	+C\sup \limits_{\theta\in S^2}\left|\int_{-\infty}^{+\infty} \int_{S^{2}}\int_{S^{2}} QTKAk^ld\theta''d\theta'dk \right|, 
	$$
	$$\label{eq:A^0}	
	\left| 	\int_{-\infty}^{+\infty} \int_{S^{2}}\int_{S^{2}} A(k,\theta',\theta)k^2dkd\theta'd\theta \right| \le \sup\limits_{x\in R^3} \left| q\right|  + C_0\left\| q\right\|_{W_2^1(R^3)}\left\| q\right\|_{L_2(R^3)}\left(  \left|\int_{S^{2}}TKAd\theta'' \right| +1  \right) 
	.$$

\end{lemma}
\begin{proof}

	Using the definition of the amplitude, Lemmas 3 and 4, and the lemma of Jordan yields 
	$$
	\int_{-\infty}^{+\infty}\int_{S^{2}}\int_{S^{2}}	A(k,\theta^{'},\theta)k^ldkd\theta'd\theta =-	\int_{-\infty}^{+\infty}\frac{1}{4\pi }\int_{S^{2}}\int_{S^{2}}\int_{R^{3}}q(x)\Psi _{+ } (k,\theta,x)e^{-ik\theta^{'} x}k^ldxdkd\theta' =	$$
	$$\label{eq:eq2}-\frac{1}{4\pi }\int_{S^{2}}\int_{S^{2}}\int_{R^{3}}q(x) \left[ e^{ik\theta x} +\sum_{n\ge  1}(-T_-D)^n\Psi_0\right]  e^{-ik\theta^{'} x}k^ld\theta' dxdk 
	$$$$
	=\int_{-\infty}^{+\infty}\int_{S^{2}}\int_{S^{2}}\tilde{q}(k(\theta-\theta'))k^ldkd\theta'd\theta+ \sum_{n\ge  1}W_n 
	,$$%	
	$$
	W_1= \int_{R^{3}}\int_{-\infty}^{+\infty}\int_{S^{2}}\int_{S^{2}}\frac{sA(s,\theta^{''},\theta)e^{-ik\theta^{'} x}q(x)e^{i s\theta'' x }}{k-s}k^l dkdxdsd\theta'd\theta'',
	$$
	
	$$
	\left| 	W_1\right| \le  C\sup \limits_{\theta\in S^2}\left|\int_{-\infty}^{+\infty} \int_{S^{2}}\int_{S^{2}} QTKAk^ld\theta''d\theta'dk \right|  
	.$$

	Similarly, 
	$$
	\left| 	W_n\right| \le C\sup \limits_{\theta\in S^2}\left|\int_{-\infty}^{+\infty} \int_{S^{2}}\int_{S^{2}} QTKAk^ld\theta''d\theta'dk \right| \left|\int_{S^{2}}TKAd\theta'' \right|^n
	.$$
	Finally,
	
	$$
	\left| 	\int_{-\infty}^{+\infty} \int_{S^{2}}\int_{S^{2}} A(k,\theta',\theta)dkd\theta'd \theta \right| \le \left| \int_{-\infty}^{+\infty}\int_{S^{2}}\int_{S^{2}}\tilde{q}(k(\theta-\theta'))dkd\theta d\theta'\right|
	$$$$
	+C_0\left\| q\right\|^2_{L_2(R^3)}\left(  \left|\int_{S^{2}}TKAd\theta'' \right| +1  \right), 
	$$
	
	$$
	\left| 	\int_{-\infty}^{+\infty} \int_{S^{2}}\int_{S^{2}} A(k,\theta',\theta)k^2dkd\theta'\right| \le \sup\limits_{x\in R^3} \left| q\right|  + C_0\left\| q\right\|^2_{L_2(R^3)}\left(  \left|\int_{S^{2}}TKAd\theta'' \right| +1  \right). 
	$$
	%%%%%%%%%%%%%%%%%%%%%%%%%%%%%%%%%%%%5

	This completes the proof.	
\end{proof}

\begin{lemma}\label{lm:l6}
	Let $$  \sup\limits_{ k}  \int_{S^{2}}\left|    \int\limits_{-\infty}^{\infty} \frac{pA(p,\theta^{'},\theta) }{ 4\pi(p-k+i0 ) }dp  \right|d\theta<\alpha<1/2,\,\,\, \sup\limits_{ k} \left| pA(p,\theta^{'},\theta)\right|<\alpha<1/2. 
	$$  Then,
	$$
	|T_{- }D\Psi _{0} |< \frac{\alpha }{1-\alpha},\,\,\,\, |T_{+ }D\Psi _{0} |< \frac{\alpha }{1-\alpha}, \,\,\,\,\,\,\,\,\,
	|D\Psi _{0} |< \frac{\alpha }{1-\alpha}, 	$$
	$$\label{eq:eq2}
	\,\, 	T_{- }g_{-}=( I-T_{- }D )^{-1 }T_{- }D\Psi _{0},
	\,\,\,\,\,\,\,\,\,\,\,\,\,\Psi _{-}=( I-T_{- }D )^{-1 }T_{- }D\Psi _{0} +\Psi _{0},           
	\label{lm:psi}
	$$
	and $q$ satisfies the following inequalities:
	$$
	\sup\limits_{x\in R^3}|q(x)|\le \left| \int_{S^{2}}TKQ d\theta\right|
	C_0\left( \left\| q\right\|^2_{L_2(R^3)}+1\right) +C_0\left\| q\right\|_{L_2(R^3)}.
	$$
\end{lemma}
\begin{proof}
	Using the equation 
	$$
	\Psi _{+ } (k,\theta,x)-\Psi _{- } (k,\theta,x)=
	%$$
	%$$
	-\frac{k}{4\pi}\int_{S^{2}}A(k,\theta^{'},\theta)\Psi_{-} (k,\theta^{'},x)d\theta^{'},~  k\in R,
	$$%
	we can write 
	$$
	T_{+}g_{+}-T_{- }g_{-}=D (T_{-}g_{-}+\Psi_{0}).
	$$
	Applying the operator $T_{-}$ to the last   equation, we have 
	$$
	T_{- }g_{-}=T_{- }D (T_{-}g_{-}+\Psi_{0}),
	$$
	$$
	(I  - T_{- }D    ) T_{- }g_{-}=T_{- }D\Psi_{0},\,\,\,
	T_{- }g_{-}=   \sum_{n\ge0}      \left( -T_{- }D \right)^{n} \Psi_{0}.
	$$
	Estimating the terms of the series, we obtain using Lemma 4
	$$
	|(T_{- }D)^{n}\Psi_{0}| \le 
	\sum_{n\ge 0 }
	\left|  \int_{-\infty}^{\infty}\dots \int_{-\infty}^{\infty} \Psi_{0}  \prod\limits_{ 0\le j<n }\frac{ \int_{S^{2}}{k_{j}}A(k_{j},\theta^{'}_{k_{j}},\theta_{k_{j}}) d\theta^{'}_{k_{j}}}{ 4\pi(k_{j+1})-k_{j}+i0 ) }   dk_{1}\dots d_{k_{n}}\right| 	$$
	$$\label{eq:eq2}\\\le \sum_{n> 0 } 2^n \alpha^{n} =\frac{2\alpha}{1-2\alpha} . 
	$$
	%	Denote $$\Lambda= \frac{1}{k^2}  \frac{\partial }{\partial k}k^2\frac{\partial }{\partial k}, r=\sqrt{x_1^2+x_2^2+x_3^2}, $$ 
	Denoting $$\Lambda=  \frac{\partial }{\partial k},\; r=\sqrt{x_1^2+x_2^2+x_3^2}, $$ 
	we have
	$$
	\Lambda\int_{S^{2}}\Psi_0d\theta =\Lambda  \frac{\sin( kr)}{ikr}=  
	\frac{\cos(kr)}{ik}-\frac{\sin( kr)}{ik^2r}, 
	$$
	$$
	\Lambda\int_{S^{2}}H_0\Psi_0d\theta =\Lambda k^2 \frac{\sin( kr)}{ikr}=  
	k\frac{\cos(kr)}{i}+\frac{\sin( kr)}{ik^2r} ,
	$$
	$$
	\left| \Lambda \int_{S^{2}}\Psi d\theta	 \right|=\left|\Lambda\int_{S^{2}}\Psi_0d\theta +\Lambda\int_{S^{2}}\sum_{n\ge 0}\left( -T_{- }D \right)^{n} \Psi_{0} d\theta\right|  > \left(\frac{1}{k} -  \frac{\alpha}{1-\alpha}\right), \;{\rm as}\,\, kr=\pi,
	$$
	and
	$$
	\Lambda\frac{1}{k-t}=
	-   \frac{1}{(k-t)^2} 
	$$

	Equation (\ref{eq:1})  yields 
	$$
	q= \frac{\Lambda\left(  H_0\int_{S^{2}}\Psi d\theta+k^2\int_{S^{2}}\Psi d\theta\right)  }{\Lambda\int_{S^{2}}\Psi d\theta} $$$$
	=\frac{2k\int_{S^{2}}T_{-}g_{-}d\theta+k^2\int_{S^{2}}\Lambda T_{-}g_{-}d\theta+  H_0\Lambda \int_{S^{2}}T_{-}g_{-}d\theta}{\Lambda\int_{S^{2}}\Psi d\theta}  $$
	$$=
	\frac{2k\int_{S^{2}}T_{-}g_{-}d\theta+\Lambda\int_{S^{2}}\sum_{n\ge 1}      \left( -T_{- }D \right)^{n}(K^2-k^2) \Psi_{0} d\theta}{\Lambda\int_{S^{2}}\Psi d\theta}$$$$
	=\frac{W_0+\sum_{n\ge 1} \int_{S^{2}}W_{n}}{\Lambda\int_{S^{2}}\Psi d\theta}.
	$$

	Denoting $$Z(k,s)=s+2k +\frac{2k^2}{k-s},$$ we then have 
	$$
	\left|  W_1\right| \le 	\left|\int_{-\infty}^{+\infty}\int_{S^{2}} \int_{S^{2}} A(s,\theta,\theta')s \frac{s^2-k^2}{(k-s)^2}  \Psi_0  \sin(\theta)ds d\theta\right|_{k=k_0} 
	$$$$
	\le \left| \int_{-\infty}^{+\infty}\int_{S^{2}}\int_{S^{2}}Z(k,)\tilde{q}(k(\theta-\theta'))\Psi_0dkd\theta\right|+
	C_0\left|\int_{S^{2}} TKQd\theta\right|.
	$$
	
	For calculating $W_{n}$, as $n \ge 1$, take the simple transformation
	$$
	\frac{s_{n}^3}{s_{n}-s_{n-1}}=	\frac{s_{n}^3-s_{n}^2s_{n-1}}{s_{n}-s_{n-1}} +\frac{s_{n}^2s_{n-1}}{s_{n}-s_{n-1}}=s_{n}^2+\frac{s_{n}^2s_{n-1}}{s_{n}-s_{n-1}}
	$$
	\begin{equation} \label{eq:6}
	=s_{n}^2+\frac{s_{n}^2s_{n-1}-s_{n}s_{n-1}^2}{s_{n}-s_{n-1}} +\frac{s_{n}s_{n-1}^2}{s_{n}-s_{n-1}}=	s_{n}^2+s_{n}s_{n-1}+\frac{s_{n}s_{n-1}^2}{s_{n}-s_{n-1}},
	\end{equation}
	$$
	\frac{As_{n}^3}{s_{n}-s_{n-1}} =As_{n}^2+As_{n}s_{n-1}+    \frac{As_{n}s^2_{n-1}}{s_{n}-s_{n-1}}=V_1+V_2+V_3 .
	$$
	Using Lemma \ref{lm:TA}   for estimating $V_1$  and $V_2$ and, for $ V_3 $, taking again the simple transformation for $s^3_{n-1} $, which will appear in the integration over  $ s_{n-1}$, we finally get 
	$$
	|q(x)|_{r=r_0}= \left| \frac{\Lambda\left(  H_0\int_{S^{2}}\Psi d\theta+k^2\int_{S^{2}}\Psi d\theta\right)  }{\Lambda\int_{S^{2}}\Psi d\theta}\right| _{k=k_0,r= \frac{\pi}{k_0}} $$$$
	\le \frac{\left| \int_{-\infty}^{+\infty}\int_{S^{2}}\int_{S^{2}}Z(k,)\tilde{q}(k(\theta-\theta'))\Psi_0dkd\theta d\theta'\right|+
		C_0\left|\int_{S^{2}} TKQd\theta\right|}{(\frac{1}{k_0}-\frac{\alpha}{(1-\alpha)} )}   + 
	$$
	$$
	% \frac{   C_0\left\| q\right\|_{W_2^1(R^3)}\left\| q\right\|_{L_2(R^3)}\left(  %\left|\int_{S^{2}}TKAd\theta'' \right| +1  \right)}{1-\frac{\alpha}{r_0^2-\frac{\alpha}{(1-\alpha)k_0}} } 
	$$	
	Finally, we get
	$$
	|q(x)|_{r=r_0}\le \sup\limits_{x\in R^3}|q(x)|\alpha+
	C_0\left\| q\right\|^2_{L_2(R^3)}+C_0\left\| q\right\|_{L_2(R^3)}+\left| \int_{S^{2}}TKQ d\theta\right|.
	$$
	
	The invariance of the Schr\"odinger equations with respect to translations and the arbitrariness of $r_0$	yield
	$$
	\sup\limits_{x\in R^3}|q(x)|\le \left| \int_{S^{2}}TKQ d\theta\right|
	C_0\left( \left\| q\right\|^2_{L_2(R^3)}+1\right) +C_0\left\| q\right\|_{L_2(R^3)}.
	$$
\end{proof}
%%%%%%%%%%%%%%%%%%%%%%%%%%%%%%%%%%%%%%%%%%%%%%%%%	

%%%%%%%%%%%%%%%%%%%%%%%%%%%%%%%%%%%%%%%%%%%%%%%%%%%%%%%%%%%%%%%%%%%%%%%%%%%%%%%%%%%%%%%%%%%%%%%%%

%------------------------------------------------
% % % % % % % % % % % % % % % % % % % % % % % % % % % % % % % % % % % % %
% % % % % % % % % % % % % % % % % % % % % % % % % % % % % % % % % % %5
% % % % % % % % % % % % % % % % % % % % % % % % % % % % % % %
% % % % % % % % % % % % % % % % % % % % % % % % % % %
\section{Discussion of the three-dimensional inverse scattering problem }
This study has shown, once again, the outstanding properties of the scattering operator, which, in combination with the analytical properties of the wave function, allows us to obtain almost-explicit formulas for the potential  from the scattering amplitude. Furthermore, this appro. The estimations following from this overcome the problem of overdetermination, resulting from the fact that the potential is a function of three variables, whereas the amplitude is a function of five variables. We have shown that it is sufficient to average the scattering amplitude to eliminate the two extra variables.
%------------------------------------------------
\section{  Studying the properties of solutions of the Cauchy problem for the Navier--Stokes equations using analytic functions generated by the Schr\"odinger equations and related to the Poincar\'e-–-Riemann–-Hilbert problem}
Numerous studies of the Navier--Stokes equations have been devoted to the problem of the smoothness of its solutions. A good overview of these studies is given in Refs. [13--17]. The spatial differentiability of the solutions is an important factor, as it controls their evolution. Obviously, differentiable solutions do not provide an effective description of turbulence. Nevertheless, the global solvability and differentiability of the solutions have not been proven, and therefore the problem of describing turbulence remains open. It is interesting to study the properties of the Fourier transform of solutions of the Navier--Stokes equations. Of particular interest is how they can be used in the description of turbulence and whether they are differentiable. The differentiability of such Fourier transforms appears to be related to the appearance or disappearance of resonance, as this implies the absence of large energy flows from small to large harmonics, which in turn precludes the appearance of turbulence.
Therefore, obtaining uniform global estimations of the Fourier transform of solutions of the Navier--Stokes equations means that the principle modelling of complex flows and related calculations will be based on the Fourier transform method. We are continuing to research these issues in relation to a numerical weather prediction model; this paper provides a theoretical justification for this approach. 

Consider the Cauchy problem for the Navier--Stokes equations:
\begin{equation}
\label{1}
%q_{t}-\nu \Delta q+\sum\limits_{k=1}^{3}q_{k}q_{x_{k}}=-\nabla
\frac{\partial\vec{v}}{\partial t}-\nu \Delta \vec{v}+(\vec{v},\nabla \vec{v})=-\nabla
p+\vec{f}(x,t),~{\rm div}~\vec{v}=0,  
\end{equation}
\begin{equation}\label{2}
\vec{v}|_{t=0}=\vec{v}_{0}(x)  
\end{equation}
in the domain $Q_{T}=R^{3}\times (0,T)$, where 
\begin{equation}
{\rm div}\;\vec{v}_{0}=0.  \label{3}
\end{equation}
The problem defined by  (\ref{1})--(\ref{3}) has at least one weak solution $ (\vec{v}, p) $ in the so-called Leray--Hopf class [16].
The following results have been  proved  [15]:
\begin{theorem}
	\textbf{} If
	$$
	\vec{v}_{0}\in W_{2}^{1}(R^{3}),\vec{f}(x,t)\in L_{2}(Q_{T}),
	$$%
	there is a single generalised solution of (\ref{1})--(\ref{3}) in the domain $Q_{T_{1}}$, $T_{1}\in \lbrack 0,T]$, satisfying the following conditions:
	$$
	\vec{v},\nabla ^{2}\vec{v},\ \ \ \nabla p\in L_{2}(Q_{T}).
	$$
	\label{thm1}
\end{theorem}
Note that $T_{1}$ depends on $\vec{v}_{0}$ and $\vec{f}(x,t)$.
\begin{lemma}\label{lm:l8}
	If we let $\vec{v_0} \in W_{2}^{2}(R^{3}),\vec{f} \in L_2(Q_T)  $,
	then the solution of (\ref{1})--(\ref{3}) satisfies  the following inequalities:
	$$
	\sup\limits_{0\leq t\leq
		T}||\vec{v}||_{L_{2}(R^{3})}^{2}+\nu\int\limits_{0}^{t}||\nabla \vec{v}||_{L_{2}(R^{3})}^{2}d\tau  \leq \ ||\vec{v}_{0}||_{L_{2}(R^{3})}^{2}+||\vec{f}||_{L_{2}(Q_{T})},	$$
	$$\label{eq:eq2}
	\sup\limits_{0\leq t\leq
		T}||\vec{\nabla v}||_{L_{2}(R^{3})}^{2}+\nu\int\limits_{0}^{t}||H_{0} \vec{v}||_{L_{2}(R^{3})}^{2}d\tau  	$$
	$$\label{eq:eq2}\leq || \nabla\vec{v}_{0}||_{L_{2}(R^{3})}^{2}+||\vec{f}||_{L_{2}(Q_{T})} +\int_{0}^{t}||(\vec{v},\nabla \vec{v})||_{L_{2}(R^{3})}||H_{0} \vec{v}||_{L_{2}(R^{3})},	$$
	$$\label{eq:eq2}
	\nu\int\limits_{0}^{t}||H_{0} \vec{v}||_{L_{2}(R^{3})}^{2}d\tau \le C+ \frac{1}{\nu}\int_{0}^{t}||(\vec{v},\nabla \vec{v})||^{2}_{L_{2}(R^{3})}dt	.
	$$
\end{lemma}
\begin{lemma}
	Let $\vec{v_0} \in W_{2}^{2}(R^{3}),$ $\vec{\tilde{v_0}} \in W_{2}^{2}(R^{3}),$ and $\vec{f} \in L_2(Q_T) $.
	Then, the solution of (\ref{1})--(\ref{3}) satisfies the following:
	$$
	\widetilde{\vec{v}}=\widetilde{\vec{v}}_{0}+
	\int\limits_{0}^{t}e^{-\nu k^{2}|(t-\tau )}(\widetilde{[(\vec{v},\nabla )\vec{v}]}+\widetilde{\vec{F}})
	d\tau ,
	\label{eqno8.4}
	$$
	where $\vec{F}=-\nabla p+\vec{f}$.
\end{lemma}
\begin{proof}
	This follows from the definition of the Fourier transform and the theory of linear differential equations.
\end{proof}

% % % % % % % % % % % % % % % % % % % % % % % % % %55
% % % % % % % % % % % % % % % % % % % % % % % % % % %5

Let us introduce the operators $F_{k}$ and $F_{kk\prime}$ as $$F_{k}f=\int_{R^{3}} e^{i(k,x)} f(x)dx,\,\,\,F_{kk\prime}f=\int_{R^{3}} e^{i(k,x)-i(x,k^{\prime})} f(x)dx, $$
$$\vec{\tilde{v}} (k)=F_{k}\vec{v},\,\, \vec{V} (k,k^{\prime})=F_{kk\prime}\vec{v}= \int_{R^{3}} e^{i(k,x)-i(x,k^{\prime})}\vec{v}dx.$$
%%%%%%%%%%%%%%%%%%%%%%%%%%%%%%%%%%%%%%%%%%%%%%%%%%%%%%%%%%%%%%%%%%%%%%%%%%%%%%%%%%%%%%%%%%%%%%%%%\

%%%%%%%%%%%%%%%%%%%%%%%%%%%%%%%%%%%%%%%%%%%%%%%%%%%%%%%%%%%%%%%%%%%%%%%%%%%%%%%%%%%%%%%%%%%%%%%%%\
\begin{lemma}\label{lm:19}
	Let $\vec{v_0} \in W_{2}^{2}(R^{3})$, $\vec{f}\in L_2(Q_T)$, and $ \left| TKV_0\right| +\left| TKV_0\right| +\left| TK^2V_0\vec{\tilde{v_0}}\right| <C$. Then,
	the solution of (\ref{1})--(\ref{3}) in Theorem \ref{thm1} satisfies  the following inequalities:
	$$
	|\tilde{v} (k)|<C,\,\,
	\,\,
	$$
	$$\label{eq:eq2}
	|TK\tilde{v} (k)|<C_0||v||_{L_2(R^3)}+   \frac{C_0t}{\sqrt{\nu}}||\nabla v||_{L_2(R^3)}||v||_{L_2(R^3)}
	.$$
\end{lemma}
\begin{proof}
	This follows from 
	$$
	\vec{\dot v}=
	-(\vec{v}\nabla )\vec{v}+(\nu  \vec{v} + \nabla p) + F ,\\
	$$	
	$$
	\vec{\tilde{v}} = \vec{\tilde{v}}_0+
	\int_{0}^{t} e^{-\nu k^2(t-\tau)}F_{k}\left( -\ (\vec{v},\nabla )\vec{v})+	  \nabla p + 	F \right)d\tau .
	$$	
	From the last equation we have
	$$
	|\vec{v} |\le |\vec{v} _0|+C_{T}.
	$$

	Denote $$\beta=\sqrt{\nu(t-\tau)},\,\,\,a=\theta x$$
	formula  121 (23) from [11] as $n=0$: yield	
	$$
	\left| 	TK\vec{v} \right| <\left| k e^{-\beta^2k^2}\right|  + \sqrt{\pi}\beta^{-1}e^{- \frac{a^2}{8\beta^2}}D_{0}\left(\frac{a}{\sqrt{2}\beta }\right), 		$$
	$$\label{eq:eq2}
	\left| TK\vec{v} \right| \le	\left| TK\vec{v}_0\right| 	$$
	$$\label{eq:eq2}	
	+\left| TK\int_{0}^{t} e^{-\nu k^2(t-\tau)}F_{k}\left( -	 (\vec{v},\nabla )\vec{v}]+ \nabla p + 	F\right)dk\right|	$$
	$$\label{eq:eq2}	
	\le \left| TK\vec{v}_0\right| +
	\int_{0}^{t}\left| k e^{-\beta^2k^2}\right|  +	\left|  \sqrt{\pi}\beta^{-1}e^{- \frac{a^2}{8\beta^2}}D_{0}(\frac{a}{\sqrt{2}\beta }) \right|||\nabla \vec{v}||_{L_{2}(R^3)}dt 
	$$$$
	\le C_0||v||_{L_2(R^3)}+   \frac{C_0t}{\sqrt{\nu}}||\nabla v||_{L_2(R^3)}||v||_{L_2(R^3)}	.
	$$

\end{proof}
%%%%%%%%%%%%%%%%%%%%%%%%%%%%%%%%%%%%%%%%%%%%%%%%%%%%%%%%%%%%%%%%%%%%%

\begin{lemma}\label{lm:19}
	Let $\vec{v_0} \in W_{2}^{2}(R^{3})$, $\vec{f}\in L_2(Q_T)$, and $ \left| TKV_0\right| +\left| TKV_0\right| +\left| TK^2V_0\vec{\tilde{v_0}}\right|
	.$ Then,
	the solution of (\ref{1})--(\ref{3}) in Theorem \ref{thm1} satisfies  the following inequalities:
	$$
	|\vec{V}  (k,k^{\prime})|<C,\,\,
	k|\vec{V}  (k,k^{\prime})|<\frac{C}{ \sqrt{(1-\cos(\theta))}}	,\,\,	$$
	$$\label{eq:eq2}	 |T\vec{V} K|<C_0||v||_{L_2(R^3)}+   \frac{C_0t}{\sqrt{\nu(1-\cos(\theta))}}||\nabla v||_{L_2(R^3)}||v||_{L_2(R^3)}.	$$
	
\end{lemma}
\begin{proof}
	This follows from 
	$$
	\vec{\dot V}=
	-F_{kk\prime}[(\vec{v},\nabla )\vec{v}]+F_{kk\prime} (\nu\Delta \vec{v} +\nabla p) + F_{kk\prime}{F} .
	$$
	
	After the transformations, we obtain
	$$
	\vec{\dot V}=
	-F_{kk\prime}[(\vec{v}\nabla )\vec{v}]+(\nu_{k} F_{kk\prime}  \vec{v} +F_{kk\prime} \nabla p) + F_{kk\prime}{F} ,\\
	$$	
	$$
	\vec{V} = \vec{V}_0+
	\int_{0}^{t} e^{-\nu k^2(1-\cos(\theta))(t-\tau)}\left( -	F_{kk\prime}[(\vec{v},\nabla )\vec{v}]+	 F_{kk\prime} \nabla p + 	F_{kk\prime}{F} \right).
	$$	
	From the last equation, we have
	$$
	|\vec{V} |\le |\vec{V} _0|+C_0\int_{0}^{t}||\nabla v||_{L_2(R^3)}||v||_{L_2(R^3)}	d\tau.
	$$

	Denote $\beta=\sqrt{(1-\cos(\theta))(t-\tau)\nu}$, $a=(\theta-\theta')x$ 	formula  121 (23) from [11] as $n=0$: yield	
	$$
	\left| 	TK\vec{V} \right| <\left| k e^{-\beta^2k^2}\right|  + \sqrt{\pi}\beta^{-1}e^{- \frac{a^2}{8\beta^2}}D_{0}\left(\frac{a}{\sqrt{2}\beta }\right) ,		$$
	$$\label{eq:eq2}
	\left| TK\vec{V} \right| \le	\left| TK\vec{V}_0\right|	$$
	$$\label{eq:eq2}	
	+\left| TK\int_{0}^{t} e^{-\nu k^2(1-\cos(\theta))(t-\tau)}\left( -	F_{kk\prime} (\vec{v},\nabla )\vec{v}]+	 F_{kk\prime} \nabla p + 	F_{kk\prime}{F} \right)dk\right|	$$
	$$\label{eq:eq2}	
	\le \left| TK\vec{V}_0\right| +
	\int_{0}^{t}\left| k e^{-\beta^2k^2}\right|  +	\left|  \sqrt{\pi}\beta^{-1}e^{- \frac{a^2}{8\beta^2}}D_{0}\left(\frac{a}{\sqrt{2}\beta }\right) \right|||\nabla \vec{v}||_{L_{2}(R^3)}|| \vec{v}||_{L_{2}(R^3)}dt 
	$$$$				
	<C_0||v||_{L_2(R^3)}+   \frac{C_0t}{\sqrt{\nu(1-\cos(\theta)) }}||\nabla v||_{L_2(R^3)}||v||_{L_2(R^3)}.
	$$

\end{proof}
%%%%%%%%%%%%%%%%%%%%%%%%%%%%%%%%%%%%%%%%%%%%%%%%%%%%%%%%%%%%%%%%%%%%%

%%%%%%%%%%%%%%%%%%%%%%%%%%%%%%%%%%%%%%%%%%%%%%%%%%%%%%%%%%%%%%%%%%%%%

%%%%%%%%%%%%%%%%%%%%%%%%%%%%%%%%%%%%%%%%%%%%%%%%%%%%%%%%%%%%%%%%%%%%%%%%%%%%%%%%%%%%%%%%%%%%%%%%%
\begin{theorem}
	\textbf{\label{Theorem 8.6}} 
	Let $\vec{v_0} \in W_{2}^{2}(R^{3}),$ $\vec{f}\in L_2(Q_T)$, $\vec{\tilde{f}}\in W_2^{2,1}(Q_T)$, $ \left| TKV_0\right| +\left| TKV_0\right| +\left| TK^2V_0\vec{\tilde{v_0}}\right| <C$,and $\int_{0}^{\infty}||  H_0\vec{f}||_{L_{2}(R^{3})}dt<C $. Then,
	the solution of (\ref{1})--(\ref{3}) in Theorem \ref{thm1} satisfies the following inequalities:
	$$
	\sup\limits_{x\in R^{3}}||\vec{v}(x)||<C,
	$$
	$$
	\left\vert \left\vert\nabla 
	\vec{v}\right\vert \right\vert
	_{L_{2}(R^{3})}+
	\nu\int\limits_{0}^{T}\int\limits_{R^{3}}|H_{0} \vec{v}
	|^{2}dxd\tau \leq {\rm const}. 
	$$
\end{theorem}
\begin{proof}
	
	Consider the Cauchy problem for the Navier--Stokes equations:
	\begin{equation}
	\label{1a}
	%q_{t}-\nu \Delta q+\sum\limits_{k=1}^{3}q_{k}q_{x_{k}}=-\nabla
	\frac{\partial\vec{v}}{\partial t}-\nu \Delta \vec{v}+(\vec{v},\nabla \vec{v})=-\nabla
	p+\vec{f}(x,t),~{\rm div}~\vec{v}=0,  
	\end{equation}
	\begin{equation}\label{2a}
	\vec{v}|_{t=0}=\vec{v}_{0}(x)  
	\end{equation}
	in the domain $Q_{T}=R^{3}\times (0,T)$, where 
	\begin{equation}
	{\rm div}\;\vec{v}_{0}=0.  \label{3a}
	\end{equation}
	We perform the following transformations:
	$$
	\vec{u_{\epsilon}}=\epsilon\vec{v}, \;p_{\epsilon}=p\epsilon,\,\,f_{\epsilon}=f{\epsilon^2},\,\,\nu_{\epsilon}=\epsilon\nu, s=\frac{t}{\epsilon}.
	$$
	Then, 
	\begin{equation}
	\label{1a}
	%q_{t}-\nu \Delta q+\sum\limits_{k=1}^{3}q_{k}q_{x_{k}}=-\nabla
	\frac{\partial\vec{u_{\epsilon}}}{\partial s}-\nu_{\epsilon} \Delta \vec{u_{\epsilon}}+(\vec{u_{\epsilon}},\nabla \vec{u_{\epsilon}})=-\nabla_{\epsilon}
	p_{\epsilon}+\vec{f_{\epsilon}}(x,t),~{\rm div}~\vec{u_{\epsilon}}=0,  
	\end{equation}
	\begin{equation}\label{2a}
	\vec{u_{\epsilon}}|_{t=0}=\vec{u_{\epsilon}}_{0}(x)  
	\end{equation}
	in the domain $Q_{T}=R^{3}\times (0, T_{\epsilon})$, where 
	\begin{equation}
	{\rm div}\;\vec{u_{\epsilon}}|_{t=0}=0.  \label{3a}
	\end{equation}
	Let us return for convenience to the notation $v_{i}=u_{\epsilon_{i}}$,
	using the equation for each   $v_{i}=u_{\epsilon_{i}} $. This gives us
	$$
	-\Delta_{x}\Psi  +v_{i}\Psi  =k^{2}\Psi, ~k\in C.  
	\label{eq:se}
	$$
	Using Lemmas 12-15,  we get estimates for 
	$$ 
	A_{i},\;\vec{V}_{i},\; TA_{i},\;T\vec{V}_{i},\; kA_{i},k\vec{V}_{i}, \;TKA_{i},\;TK\vec{V}_{i},\;TK\tilde{v_{i}},\;TK^2V\tilde{v_{i}}. $$
	The last estimations yield the representation 	
	$$ 
	q= \frac{\Lambda\left(  H_0\int_{S^{2}}\Psi d\theta+k^2\int_{S^{2}}\Psi d\theta\right)  }{\Lambda\int_{S^{2}}\Psi d\theta}|_{r=\frac{\pi}{k_0} ,k=k_0},
	$$
	and   Lemma \ref{lm:l6}	implies	
	$$\label{eq:eq2}
	\left\vert \left\vert\nabla 
	\vec{v}\right\vert \right\vert^2
	_{L_{2}(R^{3})}+\nu_{\epsilon} \int\limits_{0}^{t}||H_{0} \vec{v}||_{L_{2}(R^{3})}^{2}d\tau \le \int_{0}^{\infty}|| ( \vec{v})||_{L_{2}(R^{3})}|| ||  H_0\vec{f}||_{L_{2}(R^{3})}d\tau +   $$ $$
	\left\vert \left\vert\nabla 
	\vec{v_{0}}\right\vert \right\vert^2
	_{L_{2}(R^{3})}
	+ \frac{C_0}{\nu_{\epsilon}}\int_{0}^{t} \left(\frac{C_1}{\nu_{\epsilon}} ||(\nabla \vec{v})||^{2}_{L_{2}(R^{3})}||( \vec{v})||^{2}_{L_{2}(R^{3})} +|| \vec{v}||^{2}_{L_{2}(R^{3})}\right) ||(\nabla \vec{v})||^{2}_{L_{2}(R^{3})}d\tau .
	$$
	Denote $$
	\alpha(s)=\frac{C_0}{\nu_{\epsilon}} \left(\frac{C_1}{\nu_{\epsilon}} ||(\nabla \vec{v})||^{2}_{L_{2}(R^{3})}||( \vec{v})||^{2}_{L_{2}(R^{3})} +|| \vec{v}||^{2}_{L_{2}(R^{3})}\right) ,
	$$	
	
	$$\int_{0}^{ \frac{T}{T\epsilon\nu}}\alpha(s) ds \le \int_{0}^{\frac{1}{\nu\epsilon} }\frac{C_0}{\nu_{\epsilon}} \left(\frac{C_1}{\nu_{\epsilon}} ||(\nabla \vec{v})||^{2}_{L_{2}(R^{3})}||( \vec{v})||^{2}_{L_{2}(R^{3})} +|| \vec{v}||^{2}_{L_{2}(R^{3})}\right)ds  
	$$
	$$
	\le \frac{C_0C_1}{\nu^3_{\epsilon}}\sup\limits_{ t} || ( \vec{v})||^{2}_{L_{2}(R^{3})}\int_{0}^{\infty}\nu_{\epsilon}||(\nabla \vec{v})||^{2}_{L_{2}(R^{3})}||ds +\frac{C_0}{\nu_{\epsilon}}\sup\limits_{ t} || ( \vec{v})||^{2}_{L_{2}(R^{3})}
	$$
	$$
	\le \frac{C_0  \epsilon^4}{\epsilon\nu^3_{\epsilon}}	+
	\frac{C_0\epsilon^2 \frac{\nu}{\epsilon}}{\nu_{\epsilon}}\le 2C_0.	
	$$
	As $\epsilon=\nu \epsilon_0 $,
	the Gronwall--Bellman  lemma  yields  
	$$
	\left\vert \left\vert\nabla 
	\vec{v}\right\vert \right\vert^2
	_{L_{2}(R^{3})}+
	\nu_{\epsilon}\int\limits_{0}^{t}\int\limits_{R^{3}}|H_{0} \vec{v}
	|^{2}dxd\tau \leq  	\left\vert \left\vert\nabla 
	\vec{v_0}\right\vert \right\vert^2_ {L_{2}(R^{3})} e^{2C_0}$$$$+
	e^{2C_0}\int_{0}^{\infty}|| ( \vec{v})||_{L_{2}(R^{3})}|| ||  H_0\vec{f}||_{L_{2}(R^{3})}d\tau .
	$$
\end{proof}
%%%%%%%%%%%%%%%%%%%%%%%%%%%%%%%%%%%%%%%%%%%%%%%%%%%%%%%%%%%%%%%%%%%%%
Theorem \ref{Theorem 8.6}  \,\,\,asserts the global solvability and uniqueness of the Cauchy problem for the Navier--Stokes equations.

\section{Discussion}
As noted in the introduction, the key method of investigating the Cauchy problem for the Navier--Stokes equations is its reduction to the Poincar\'e--Riemann--Hilbert problem. By studying the wave functions for the Schr\"dinger equation of the generated velocity components, we obtain unique estimates for the maximum velocity.
Uniform global estimations of the Fourier transform of solutions of the Navier--Stokes equations indicate that the principle modelling of complex flows and related calculations can be based on the Fourier transform method. In terms of the Fourier transform, under both smooth initial conditions and right-hand sides, no exacerbations appear in the speed and pressure modes. A loss of smoothness in terms of the Fourier transform can only be expected in the case of singular initial conditions or of unlimited forces in $ L_ {2} (Q_ {T}) $.
The theory developed by us is supported by numerical calculations performed in Refs. [18--20],
where the dependence of the smoothness of the solution on the oscillations of the system is clearly deduced.

\section {Reduction of the Riemann hypothesis to the Poincar\'e--Riemann--Hilbert problem}
%\section{The Elsevier article class}
This study is concerned with the properties of modified zeta functions. Riemann's zeta function is defined by
the Dirichlet series
%$$
%\section{The Elsevier article class}
This study is concerned with the properties of modified zeta functions. Riemann's zeta function is defined by
the Dirichlet series
%$$
\begin{equation}
\zeta(s)=\sum_{n=1}^{\infty}\frac{1}{n^s},\,\ s=\sigma +it,
\end{equation}
%$$
which is absolutely and uniformly convergent in any finite region of the complex $s$-plane for which
$ \sigma \ge 1 + \epsilon,\epsilon > 0.$ If $\sigma > 1$, then $\zeta$ is represented by the following Euler product formula
\begin{equation}
\zeta (s) =   \prod_{p} \left[  1- \frac{1} {p^s}  \right]^{-1},  
\end{equation}
where $p$ runs over all prime numbers.
$\zeta(s)$  was first introduced by Euler in 1737 [21], who also obtained formula (2). Dirichlet and Chebyshev considered this function in their study on the distribution of prime numbers [22]. 
However, the most profound properties of $\zeta(z)$ were only discovered later, when it was extended to the complex plane. In 1876, Riemann [23] proved that
$\zeta(s)$ allows analytical continuation to the entire $z$-plane as follows:
%$$
\begin{equation}
\pi^{-s/2} \Gamma(s/2) \zeta(s) = 1/(s(s-1)) +  \int\limits_{1}^{+\infty}  (x^{s/2-1}  + x^{-(1+s)/2} )\theta(x) dx,
\end{equation}
%$$
where $ \Gamma (z) $ is the gamma function and 
$$
\theta(x) = \sum_{n=1}^{\infty} \exp (-\pi n^{2}x).
$$$\zeta(s)$ is a regular function for all values of $s$, except $s=1$, where it has a simple pole with residue $1$; 
moreover, it satisfies the following functional equation:

\begin{equation}
\pi^{ -s/2} \Gamma(s/2) \zeta(s) =\pi^{-(1-s)/2} \Gamma((1-s)/2) \zeta (1-s)
\end{equation}This equation is called Riemann's functional equation.

Riemann's zeta function is an important subject of study and has numerous interesting generalizations.
The role of the zeta function is highly significant in number theory, where it is connected with various fundamental functions, such as the M\"obius function, the Liouville function, the
number of divisors, and the number  of prime divisors. The detailed theory
of zeta functions is presented in [24]. The zeta function has found application in various other fields, notably in quantum statistical mechanics and quantum field theory [25--27]. Riemann's zeta function is often introduced in quantum statistics formulas. A well-known example is 
the Stefan-Boltzman law of a black body's radiation. Its ubiquitous use in seemingly unrelated areas demonstrates the necessity for further investigation.

The present study is concerned with the analytical properties of the
following generalized zeta functions:

$$
P (s)=\sum_{j\ge 1}\frac{1}{p_{j}^s} ,\,Re(s)>1+\delta, \delta>0,
$$

where $ \{p_ {j}:j\geq1\} $ is an increasing enumeration of all prime numbers. The form of  $P (s)$  suggests that it possesses
the same properties as the zeta function; however, this is not quite obvious and can be seen by considering.
$$
ln ( \zeta (s))  = \sum_{n=1}^\infty P(ns)/n,\,  \,Re(s)>1+\delta, \delta>0.
$$

%%%%%%%%%%%%%%%%%%%%%%%%%%%%%%%%%%%%%%%%%%%%%%%%%%%%%%%%%%%%%%%%%%%%%%%%

Hadamard was the first to apply $P (s)$ in the study of the zeta function [28]. Chernoff made significant progress in the Riemann hypothesis using $P (s)$ [29]. In the present study, modifications of Chernoff's results are obtained. Specifically, his study on the pseudo zeta function is completed. Chernoff obtained an equivalent formulation of the Riemann hypothesis in terms of a pseudo zeta function as follows.

\bigskip
\noindent
THEOREM.  (Chernoff)
Let $$C(s) =\prod_{n>1} \left[  1- \frac{1} {(nln(n))^s}  \right]^{-1}$$.
Then, $C(s)$ continues analytically into the critical strip and has no zeros there.

\bigskip
The significance of this theorem is that if the primes were distributed
more regularly (i.e., if  $p_{n} \equiv n \log n $), then the Riemann
hypothesis would be trivially true.
In an effort to further develop the work of Chernoff and Hadamard, the following question naturally arises: Does the pseudo zeta function $P(s)$   continues analytically into the critical strip? 
It should be noted that analytic extensions of $P (s)$ were first studied by E. Landau and A. Walvis [30] and T. Estarmann [31], [32];	 however, no satisfactory estimates for $P(s)$ were obtained, and the present study is concerned with this question.
\newpage
THEOREM.  ( E. Landau,\,A. Walvis,\,T. Estarmann) \\ Let
$ \mu(n)-function M\ddot{o}bius \,\,\,$.
Then,$$
P(s)=\sum_{n\ge 1 }    \frac{\mu(n)ln \zeta(ns)}{n} \,\,as \,\,\, Re(s)>1+\delta,\,\, \delta>0,
$$		$$
P_0(s)=	\sum_{n\ge 1 }    \frac{\mu(n)ln\zeta(ns)}{n}-meromorphic\,\,function\,\,as\,\, Re(s)>\delta,\,\, \delta>0.	$$\\
We introduce the following analogs of the function 	P(s):

$$
Q_2 (s)  = ln ( \zeta (s))-\sum_{n=m}^\infty P(ns)/n,\,   \,Re(s)>1/2+\delta,\\
$$
\begin{equation}	
Q_2 (1-s)  = ln ( \zeta (1-s))-\sum_{n=m}^\infty P(n(1-s))/n,\,   \,Re(s)<1/2-\delta\\
\end{equation}

	The paper is organized as follows. Intermediate estimates are first  obtained for the $ ln\zeta(s) $. Subsequently, the sets where the logarithm of the zeta function is uniquely determined are defined. These sets are composed of rectangles in which the zeta function has no roots, and they cover the entire critical strip except for the rectangular regions in which the zeros of the zeta functions are located. In the rectangles in which there are no zeros of the zeta function, the real value of its logarithm can be defined, and in these sets, the mirror-symmetric equation that arises by taking the logarithm on both sides of the Riemann functional equation is investigated. Then, the Fourier transform is applied to it, and it is multiplied by a regulating factor. Thus, a Riemann-Hilbert boundary value problem is obtained for the $ Q_{2}(s) $
The properties of the solution to the Riemann--Hilbert boundary value problem are expressed in terms of the Hilbert integral transform.
In the rectangles in which the zeta function has no roots, the Hilbert transform can be used to obtain exact lower bounds for the zeta function in the critical strip.

\section{RESULTS}
As mentioned in Introduction, certain simple intermediate estimates are first obtained.

%	Subsequently, two types of rectangles are constructed: those in which the zeta function has no zeros, and those in which it has one zero.
%	This step is very important, as these rectangles will be repeatedly used in the sequel. The goal is to compare the values ​​of the zeta function in adjacent rectangles.
The rectangles in which the zeta function hasn't zeros are first introduced as follows:
%	The rectangles  are first introduced as follows:
$$
\,\,D_+(n,\epsilon)=(  s|1/2+\epsilon< Re(s)<3/2-\epsilon,\, Im(s_{n}) <Im(s)<Im(s_{n}) +d_n  , 1<|Im(s)| )
$$
$$
D_-(n,\epsilon)=(  s|1/2+\epsilon< Re(s)<3/2-\epsilon,\, -Im(s_{n}) -d_n <Im(s) <-Im(s_{n})  , 1<|Im(s)|) 
$$
Where 
$$\zeta(s_{n+1})=0,\,\zeta(s_n)=0,\, \zeta(1-s_n)=0,\zeta(1-s_{n+1})=0,\,\,\zeta(1-s_n)=0, \\
$$
$$
d_n=(Im(s_{n+1})-Im(s_{n})),\,\,\,|s_n-s_{n+1}|>0	$$
The condition $ 1 <| Im (s) | $ is necessary to exclude the point of the pole of the zeta function. This restriction does not limit our reasoning, as it was shown in [35] that the zeros of the zeta function closest to the real axis lie on the critical line.

The sets of $D(n,\epsilon),$ is shown in the figure below.
\\
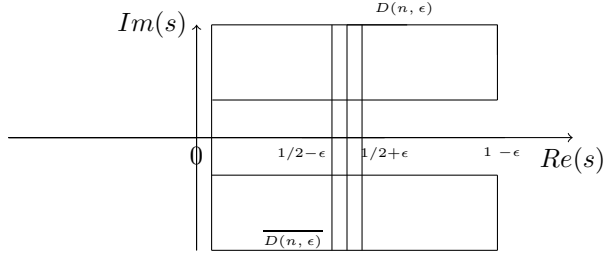
\begin{tikzpicture}
\draw (-2.5,0) -- (2.5,0);
\draw (2*1.0,1.5) -- (2*2,1.5);
\draw[->] (-2.5,0) -- (2*2.5,0) node[anchor=north] {$Re(s)$};
\draw[->] (0,-1.5) -- (0,1.5) node[anchor=east] {$Im(s)$};
\draw[] (-2.5,0) -- (0,0) node[anchor=north] {$0$};
%\draw[] (-2.5,0) -- (2*1.3,0) node[anchor=north] {$ $\tiny {1/2}$$};
%\draw[] (-2.5,0) -- (2*1.5,0) node[anchor=north] {$ $\tiny {+}$$};

\draw[] (2.2,0) -- (2*1.25,0) node[anchor=north] {$$\tiny {1/2}$+ \epsilon $};
%(-2.5,1.5)  (2*1.2,1.5) {$  \epsilon $};
%\draw[] (-2.5,0) -- (2*0.9,0) node[anchor=north] {$ $\tiny {1/2}$$};
\draw[] (2.2,0) -- (2*0.7,0) node[anchor=north] {$$\tiny {1/2}$- \epsilon $};
%\draw[] (-2.5,0) -- (2*1.8,0) node[anchor=north] {$ $\tiny {1}$$};
%\draw[] (-2.5,0) -- (2*1.75,0) node[anchor=north] {$  \Omega_n_{\epsilon,\epsilon}$};
\draw[] (-2.5,0) -- (2*2.05,0) node[anchor=north] {$ $\tiny {1 $- \epsilon $}  $  $};

%\draw (2*1.1,-1.5) -- (2*1.9,-1.5);
\draw (0,-1.5) -- (0,1.5);
\draw (2*1,-1.5) -- (2*1,1.5);
\draw (2*1.1,-1.5) -- (2*1.1,1.5);
\draw (2*2,1.1) -- (2*2,1.5);
%	\draw (2*0.1,1.1) -- (2*1.4,1.1)
\draw (2*0.1,1.5) -- (2*1.4,1.5)node[anchor=south] {$ $\tiny {$D(n,\epsilon)$}  $  $};
%\draw (2,1.1) -- (4,1.1);
%\draw[] (2*1.2,0) -- (0,0) node[anchor=north] {$0$};
%\draw (2*2,0.5) -- (2*2,1);
\draw (2*0.1,-1.5) -- (2*0.1,1.5);
\draw (2*0.9,-1.5) -- (2*0.9,1.5);
%%%%%%%%%%%%%%%%%%%%%%%%%%%%
%	\draw (2.8,1) -- (2.811,1)node[anchor=north] {$ $\tiny {$D_0(n,\delta)$}  $  $};

%	\draw (1,-1) -- (1.11,-1)node[anchor=south] {$ $\tiny {$\overline{D_0(n,\delta)}$}  $  $};

\draw (4,0.5) -- (4,1.5);
\draw (0.21,0.5) -- (4,0.5);
\draw (0.2,-0.5) -- (4,-0.5);
\draw (1.2,-1.1)  (1.32,-1.1)node[anchor=north]{$ $\tiny {$\overline{D(n,\epsilon)}$}  $  $};
%\draw (0.,-1.1) -- (2,-1.1);

\draw (0.2,-1.5) -- (4,-1.5);
\draw (4,-0.5) -- (4,-1.5);

%\draw[] (-2.5,0) -- (2*0.29,0) node[anchor=north] {$  \overline{\Omega_n}_{\epsilon,\epsilon}$};
\end{tikzpicture}\\
%	\newpage
\begin  {theorem}
Let $ s\in Z $  ,\,\,\,	$
F(s)= \frac{s}{2}ln(\pi)-ln(\Gamma(s/2))-  \frac{1-s}{2}ln(\pi)+ln(\Gamma(1-s)/2)).
$
Then, 

\begin{eqnarray*}
	\sup_{s\in D_+(n,\epsilon)\cup{D_-(n,\epsilon)}}|F(\tau+i\alpha)| +\sup_{s\in D_+(n,\epsilon)\cup{D_-(n,\epsilon)}}| \frac{dF(\tau+i\alpha)}{d\tau}| 
	<CC_{n}
\end{eqnarray*}
\end {theorem}

\begin{proof}
	As $ \epsilon<\tau<1-\epsilon$   implies that
	F - is holomorphic 	which completes the proof.	
\end{proof}

As mentioned in Introduction, a Riemann--Hilbert boundary value problem should be obtained. To this end, an equation should be derived that determines the difference between the boundary values of the analytic functions in the upper plane and the lower plane.

\begin{definition}	 
	let's call the regular solution of the functional equation (19) the analytic function $ \ nu (s) $ in the domain $ \, \, D _ + (n, \ epsilon) $ satisfies the following condition:\\
	1. $ ln \nu (s) = \sum_ {m = 1} ^ \infty | P (mz) / m, $ $ Re (s)> 1+ \epsilon $ \\
	2. $
	\pi ^ {-s / 2} \Gamma (s / 2) \nu (s) = \pi ^ {- (1-s) / 2} \Gamma ((1-s) / 2) \nu (1 -s)
	$
	
\end{definition}
$\nu(s)=const  \zeta(s)$ - example of irregular solution of the equation( 4).\\
\begin  {lemma}	
In domain 	$
\,\,D_+(n,\epsilon)
$
there is a unique regular solution to the equation (19) . 
\end{lemma}
\begin{proof}
The existence follows from the fact that the zeta function satisfies this equation.
Consider the difference between the two solutions, according to the conditions of the Lemma, for the difference we have an analytic function which is zeroable on a set of nonzero measure, from which its identity is zero.\end{proof}
\begin  {lemma}	

Let  $Re(s)>1/2+\epsilon$ 
then
\begin{eqnarray*}
\sum_{m=2}^\infty| P(ms)/m | <  CC_{\epsilon} \,\\
\end{eqnarray*}
\end{lemma}
\begin{proof}
The estimates of the harmonic series give the following estimates

$$
\sum_{m=2}^\infty| P(mz)/m |  \le \sum_{m=2}^\infty| P(mz)/m | \le C_{\epsilon}\sum_{m=2}^\infty|-2^{m\epsilon} /m |< CC_{\epsilon}
$$
\end{proof}

As mentioned in the introduction, one should obtain the Riemann – Gilbert boundary value problem. To do this, an equation must be derived that determines the difference between the boundary values of the analytical functions in the upper plane and the lower plane.

\begin  {theorem}
as
$$
s \in  D_+(n,\epsilon),\,\,	F_2(s)= Re\left(  F(s)-\sum_{n=2}^\infty P(ns)/n \right)  $$  

implies
\begin{eqnarray*}
Q_{2}(s) =  ln  |\zeta(1 - s)| + F_2( s),\\
\sup_{s\in D_+(n,\epsilon)}	|F_2(s)|<C_2C_n.
\end{eqnarray*}

\end {theorem}

\begin{proof}

As $ s\in D _ + (n, \ epsilon) $ (20) yield the equation for $ Q_ {2} (s) $ . The estimate for $ F_ {2} (s) $ follows from Lemma 17.
\end{proof}	
For the formulation of the Riemann – Hilbert boundary value problem, it is necessary to perform a number of transformations. Some preliminary arguments suggest that the Fourier transform is an appropriate choice. In the Riemann – Hilbert boundary value problems, the asymptotic behavior of unknown functions is very important. To ensure this behavior, it is necessary to have an assessment.
We introduce the functions $ R (k) $, $ Q_ \epsilon (s) $, $ L_\epsilon (1-s) $, $ F_\epsilon (s) $ and their Fourier transform. Below we will also use the Heaviside function.  $ \theta(x)  $ , $( \theta(x)=0, x<0 ; \theta(x)=1, x\ge 0 )  $

\begin{eqnarray*}
Q_\epsilon(s)= Q_{2}(s)\theta(Re(s)-1/2-\epsilon)),\,\, \,\,   L_\epsilon(1-s)= ln\left|  \zeta(1-s)\right| \theta(Re(s)-1/2-\epsilon)\,\, \\ R(k) = \frac{e^{-ik}}{k-ia}+1     \,\,,F_\epsilon(s)=F_2(s)\theta(Re(s)-1/2-\epsilon);
\end{eqnarray*} 
\begin{eqnarray*}
J_{\epsilon}(k,\alpha )= \frac{1}{\sqrt{2\pi}}\int_{\epsilon}^{1-\epsilon} L_\epsilon(\tau-i\alpha )e^{ik\tau}d\tau\,\,\,,
I_{\epsilon}(k,\alpha )=\frac{1}{\sqrt{2\pi}}\int_{\epsilon}^{3/2-\epsilon} Q_\epsilon(\tau+i\alpha )e^{-ik\tau}d\tau.
\end{eqnarray*}
\begin{eqnarray*}
\frac{1}{\sqrt{2\pi}}\int_{\epsilon}^{3/2-\epsilon} L_\epsilon(1-\tau-i\alpha )e^{-ik\tau}d\tau = \frac{e^{-ik}}{\sqrt{2\pi}}\int_{\epsilon}^{1-\epsilon} 
L_\epsilon(\tau-i\alpha )e^{ik\tau}d\tau+\\ {\sqrt{2\pi}}\int_{1-\epsilon}^{3/2-\epsilon} 
L_\epsilon(1-\tau-i\alpha )e^{-ik\tau}d\tau=	e^{-ik\tau}J_{\epsilon}(k,\alpha )+ 	S_{\epsilon}(k,\alpha ) ;
\\		
\widetilde{Q_\epsilon}(k,\alpha ) =\frac{1}{\sqrt{2\pi}}\int_{\epsilon}^{3/2-\epsilon} Q_\epsilon(\tau+i\alpha )e^{-ik\tau}d\tau,\,\widetilde{F_\epsilon}(k,\alpha) =\frac{1}{\sqrt{2\pi}}\int_{\epsilon}^{3/2-\epsilon}F_\epsilon(\tau +i\alpha)e^{-ik\tau}ds.
\end{eqnarray*}

To obtain the Riemann-Hilbert boundary value   problem , the following lemma is required.

\begin{lemma}
Let $ a>2$\,\,\, then  $ind (R)=0$ \,
\end{lemma}

\begin{proof}
By definition 
$$
ind(R)= \frac{1}{2\pi i}  \int_{-\infty}^{+\infty} \frac{R(k)^\prime }{R(k)}dk
$$
As $$ Im(k)<0,	|e^{-ik}|<1 \,\,and \,\,|k-ia|>2\,\, yield \,\,\, \frac{R(k)^\prime }{R(k)}$$
have nothing pole. Latest statement and Lemma of Jordan yield $ind(R)=0$.

\end{proof}

To obtain the necessary asymptotics, the following lemma is required.
\begin{lemma}
Let 
$$ a>2  $$ then
$
ln(R(k))$ is single-valued analitical fuction in lower half plane.\\
\end{lemma}
\begin{proof}
As $$ Im(k)\le 0\,\,\ \, $$ yeild 
$$Re(R(k))=1+ Re\left[ \frac{e^{ik}}{k-ia}\right]  >0$$

which completes proof.
\end{proof}
Denote $\Omega_n(s)={\rm Re}(s-s_{n})   $ and $\omega(s)=\frac{(s-s_{n})^{\rho(s_{n})}(s-1+s_{n})^{\rho(s_{n})} . }{(1-s)}  $

$ \rho(s_{n})$	is multiplicity  root of $\zeta(s)$ as $ s=s_n$
The following presents results of  [34]\\
\textbf {Theorem of Backlund R.} \\
Let $ \zeta(s_n)=0$ 
then $$\rho(s_{n})<C_0ln|s_n|. $$

All the arguments given below are based on the assumption of the error of the Riemann hypothesis i.e $ |\Omega_n(1/2)|=\epsilon_{n}>0$. At the end of our work we will encounter a contradiction with our assumption from which the truth of the Riemann hypothesis will follow
\begin{lemma}
Let $ \gamma_n = \frac{1}{4\rho(s_{n}) }  $, $ |\Omega_n(1/2)|=\epsilon_{n}>0$, and	
$d_n=({\rm Im}(s_{n+1})-{\rm Im}(s_{n}))/2	  .$ \\

Then, we have the following estimate as $\epsilon=0.01\epsilon_{n}(1-Re(s_n))$:
\begin{eqnarray*}
\sup\limits_{ Im s_{n}<\alpha<Im{ s_{n} +d_n}}  \int_{\epsilon}^{3/2-\epsilon}\left| Q_\epsilon (\tau+i\alpha)\right|^2 +\left| Q_\epsilon (\tau+i\alpha)\right|d\tau  <  C_nC_{\epsilon_{n}} C_{\gamma_n}. 
\end{eqnarray*}
\begin{eqnarray*}
\sup\limits_{ -Im s_{n}<\alpha<-Im{ s_{n} +d_n}}   \int_{\epsilon}^{3/2-\epsilon}\left| L_\epsilon (\tau-i\alpha)\right|^2 +\left| L_\epsilon (\tau-i\alpha)\right|d\tau  <  C_nC_{\epsilon_{n}} C_{\gamma_n}. 
\end{eqnarray*}

\begin{eqnarray*}
\sup\limits_{ Im s_{n}<\alpha<Im{ s_{n} +d_n}}  \int_{\epsilon}^{3/2-\epsilon}\left| F_\epsilon (\tau+i\alpha)\right|^2+\left| F_\epsilon (\tau+i\alpha)\right|d\tau <  C_nC_{\epsilon_{n}} C_{\gamma_n}. 
\end{eqnarray*}
\begin{eqnarray*}
\sup\limits_{ -Im s_{n}<\alpha<-Im{ s_{n} +d_n}}  \int_{1-\epsilon}^{3/2-\epsilon}\left| S_\epsilon (\tau-i\alpha)\right|^2+\left| S_\epsilon (\tau-i\alpha)\right|d\tau <  C_nC_{\epsilon_{n}} C_{\gamma_n}. 
\end{eqnarray*}

\end{lemma}

\begin{proof}

by definion $ Q_\epsilon (s) $, 
$$
I_Q=\int_{\epsilon}^{3/2-\epsilon}\left| Q_\epsilon (\tau+i\alpha)\right|^2 +\left| Q_\epsilon (\tau+i\alpha)\right|d\tau <$$
$$
\int_{\epsilon}^{3/2-\epsilon}	\theta(Re(s)-1/2-\epsilon)\left( \left| ln| \zeta(\tau+i\alpha)| \right|^2 +\left|ln\left|   \zeta(\tau+i\alpha)\right|  \right|+
|\sum_{n=2}^\infty P(ns)/n| \right) d\tau  \le
$$
$$	C_\epsilon+\int_{\epsilon}^{3/2-\epsilon}\left| ln \left|  \frac { \zeta(\tau+i\alpha) }{\omega(\tau+i\alpha)}\right|   \right|^2
+\left| ln \left|  \frac { 1}{\omega(\tau+i\alpha)}\right|   \right|^2 +\left| ln \left|  \frac { \zeta(\tau+i\alpha) }{\omega(\tau+i\alpha)}\right|   \right|
+\left| ln \left|  \frac { 1}{\omega(\tau+i\alpha)}\right|   \right|d\tau
$$
Denote
%	$$
%	I_d= \left| \frac{d\zeta}{dz} \right|_{z=s_{n-1}}+ 
%	 \left| \frac{d\zeta}{dz}\right| _{z=s_n}+
%	  \left| \frac{d\zeta}{dz}\right| _{z=\overline{s_{n-1}}}+  %\left|\frac{d\zeta}{dz}\right|_{z=\overline{s_n}}
%		$$

$$
L_{max}=\max\limits_{s\in D_+(n,\epsilon)\cup D_-(n,\epsilon)}\left| \frac{\zeta(s) }{\omega(s)}\right| ,\,\,\,	L_{min}=\min\limits_{s\in D_+(n,\epsilon)\cup D_-(n,\epsilon)}\left| \frac{\zeta(s) }{\omega(s)}\right| 
$$

$$
I_Q<C_\epsilon+\left|ln \left| L_{max}+\frac{1}{L_{min}}\right|  \right|+\left|ln \left| L_{max}+\frac{1}{L_{min}}\right|  \right|^2 +C_{\gamma_n}\int_{\epsilon}^{3/2-\epsilon}\left|   \frac { 1}{\omega(s)}  \right|^{2\gamma_n}
+\left|   \frac { 1}{\omega(s)} \right|^{\gamma_n} d\tau
$$
which completes the proof.
\end{proof}

%\begin{Discussion)
%\end{Discussion)
The previous constructions allow the calculation of the asymptotics as follows.
%%%%%%%%%%%%%%%%%%%%%%%%%%%%%%%%%%%%%%%%

%%%%%%%%%%%%%%%%%%%%%%%%%%%%%%%%%%%%%%%%%%%%%%%
\begin{lemma}\label{lm:l1}
Let $ (3/4+i\alpha )\in D(n,\epsilon) $,$\Omega_n(1/2)=\epsilon>2/m ,  $.Then \\
\begin{eqnarray*}
\lim_{Im(k)\to -\infty} I_\epsilon(k,\alpha ) = 0 , \,\,   \lim_{Im(k)\to \infty} J_\epsilon(k,\alpha ) = 0.
\end{eqnarray*}
and as Im(k)=0 
\begin{eqnarray*}
\lim_{Re(k)\to \infty} I_\epsilon(k,\alpha ) = 0 , \,\,   \lim_{Re(k)\to \infty} J_\epsilon(k,\alpha ) = 0.
\end{eqnarray*}

\end{lemma}	

\begin{proof}
To study the asymptotics, by Lemma 5 and the finiteness of $ \mu_\epsilon $ yield  

\begin{eqnarray*} 
\left| I_\epsilon(k,\alpha )\right| =\left| \int_{\epsilon}^{3/2-\epsilon} Q_\epsilon(\tau+i\alpha )e^{-ik\tau}d\tau\right|\le
C\int_{\epsilon}^{3/2-\epsilon} \left(  \left|   Q_\epsilon(\tau+i\alpha )\right|^2d\tau\right)^{1/2} \frac{1}{|Im(k)|^{1/2}} %=\frac{1}{(ik)^3}\int_{0}^{1}  \frac{d^3  Q_\epsilon(\tau+i\alpha )}{d\tau ^3}e^{-ik\tau}d\tau.
\end{eqnarray*}

A similar argument is used for the function
$$
J_{\epsilon}(k,\alpha )=	\frac{1}{\sqrt{2\pi}}\int_{\epsilon}^{3/2-\epsilon} Q_\epsilon(\tau-i\alpha )e^{ik\tau}d\tau. 
$$
As 	$Im(k)>0\,\,, J_\epsilon(\tau,\alpha )$ can be estimated using the last expression and Lemma 5 as follows:
$$
\left| J_\epsilon(k,\alpha ) \right|<	\int_{\epsilon}^{3/2-\epsilon} \left(  \left|  Q_\epsilon(\tau-i\alpha)\right|^2d\tau\right)^{1/2} \frac{1}{|Im(k)|^{1/2}}
$$

As Im(k)=0, by  the Riemann-Lebesgue lemma yield
\begin{eqnarray*}
\lim_{Re(k)\to \infty} I_\epsilon(k,\alpha ) = 0 , \,\,   \lim_{Re(k)\to \infty} J_\epsilon(k,\alpha ) = 0.
\end{eqnarray*}

which completes the proof.
\end{proof}

The reduction to a Riemann--Hilbert boundary value problem can now be formulated as follows.

\begin  {theorem}
Let \begin{eqnarray*} (3/4+i\alpha) \in D(n,\epsilon),a>2,\Omega_n(1/2)=\epsilon>2/m  \\
\Gamma_{+}(k)=	-\frac{1}{2\pi i}  \int_{-\infty}^{ \,\,\infty}\frac{ln (R(t))dt}{ t-k-i0}  \\
\Gamma_{-}(k)=	-\frac{1}{2\pi i}  \int_{-\infty}^{ \,\,\infty}\frac{ln (R(t))dt}{ t-k+i0}  \\
X_+(k)=e^{\Gamma_+(k)},\,\,\,
X_-(k)=e^{\Gamma_-(k)} ,\,\, R(k)=\frac{X_-(k)}{X_+(k)},\, G_{\epsilon}(k,\alpha )=J_{\epsilon}(k,\alpha ).
\end{eqnarray*}
Then,

$$ 
J_\epsilon(k,\alpha ) = -\frac{X_+(k)}{2\pi i}  \int_{-\infty}^{ \,\,\infty} \frac{G_\epsilon(t,\alpha)}{ X_-(t)}\frac{dt}{ t-k-i0} =X_+(k)T_+ \frac{G_\epsilon}{X_-}    
$$
$$ 
\frac{I_\epsilon(k,\alpha )}{k-ia}	-\frac{\widetilde{F_\epsilon}(k,\alpha)}{k-ia}	 =  -\frac{X_-(k)}{2\pi i}  \int_{-\infty}^{ \,\,\infty} \frac{G_\epsilon(t,\alpha)}{ X_-(t)}\frac{dt}{ t-k+i0}    dt=X_-(k)T_- \frac{G_\epsilon}{X_-}   .
$$
%$$ |Q(s)|<C_\epsilon $$
\end {theorem}

\begin{proof}
By Theorem 7 and  Lemma 17 we have
\begin{eqnarray*}
Q_\epsilon( s) =L_\epsilon(1 - s)+F_\epsilon(s). 
\end{eqnarray*}
Using the Fourier transform, we obtain
\begin{eqnarray*}
I_\epsilon(k,\alpha )=e^{-ik}J_\epsilon(k,\alpha )  +\widetilde{F_\epsilon}(k,\alpha)+S_\epsilon(k,\alpha )  .
\end{eqnarray*}
Multiplying this equation by $ \frac{1}{k-ia}$  we get 
\begin{eqnarray*}
\frac{I_\epsilon(k,\alpha )}{k-ia}	 =\frac{e^{-ik}J_\epsilon(k,\alpha )}{k-ia}   +     \frac{\widetilde{F_\epsilon}(k,\alpha)+S_\epsilon(k,\alpha )}{k-ia}.
\end{eqnarray*}
Rewriting latest equation
\begin{eqnarray*}
\frac{I_\epsilon(k,\alpha )}{k-ia}-\frac{\widetilde{F_\epsilon}(k,\alpha)+S_\epsilon(k,\alpha )}{k-ia}	 =R(k)J_\epsilon(k,\alpha )   + J_\epsilon(k,\alpha ).
\end{eqnarray*}

\begin{eqnarray}
\Psi_{-}(k,\alpha ) =	\frac{I_\epsilon(k,\alpha )}{k-ia}-\frac{\widetilde{F_\epsilon}(k,\alpha)+S_\epsilon(k,\alpha )}{k-ia}	\\ 	\Psi_{+}(k,\alpha ) =J_\epsilon(k,\alpha ) 
\end{eqnarray}
$$		G_\epsilon(k,\alpha)= J_\epsilon(k,\alpha )  $$

Using Lemma 20, the following Riemann-Hilbert boundary value  problem is obtained regarding the definition of an analytic function from its boundary values on the real line:
\begin{eqnarray}
\Psi_{-}(k,\alpha ) = R(k) \Psi_{+}(k,\alpha ) +G_\epsilon(k,\alpha ),\,\,\\	
\lim_{Re(k)\to \infty} 	\Psi_{+}(k,\alpha )= 0  \,\,\,as \,\,\,Im(k)\ge 0, \,\,   \lim_{Re(k)\to -\infty}   \Psi_{-}(k,\alpha )= 0 \,\,\,as\,Im(k)\le 0
\end{eqnarray}

Hilbert's formula and Lemma 19-Lemma 21 gives the solution to the Riemann-Hilbert boundary value  problem (23),(24) 
\begin{eqnarray}
\Psi_{+}(k,\alpha )=  -\frac{X_+(k)}{2\pi i}  \int_{-\infty}^{ \,\,\infty} \frac{G_\epsilon(t,\alpha)}{ X_-(t)}\frac{dt}{ t-k-i0}        \\
\Psi_{-}(k,\alpha )=  -\frac{X_-(k)}{2\pi i}  \int_{-\infty}^{ \,\,\infty} \frac{G_\epsilon(t,\alpha)}{ X_-(t)}\frac{dt}{ t-k+i0}         
\end{eqnarray}
Denote 

%%%%%%%%%%%%%%%%%%%%%%%%%%%%%%%%%%%%%%
\begin{eqnarray*}
\Phi_{+}(k,\alpha )=\Psi_{+}(k,\alpha )  -	J_\epsilon(k,\alpha )    \\
\Phi_{-}(k,\alpha )= \Psi_{-}(k,\alpha )	 -	 \frac{I_\epsilon(k,\alpha ) -\widetilde{F_\epsilon}(k,\alpha)-S_\epsilon(k,\alpha ) }{ (k-ia)} 	
\end{eqnarray*}

Considering the difference between the two solutions (23) - (24) we obtain the Riemann-Hilbert boundary value problem:
\begin{eqnarray*}
\Phi_{-}(k,\alpha ) = R(k) \Phi_{+}(k,\alpha ) \\	
\lim_{Re(k)\to \infty} 	\Phi_{+}(k,\alpha )= 0  \,\,\,as \,\,\,Im(k)>0, \,\,   \lim_{Re(k)\to -\infty}   \Phi_{-}(k,\alpha )= 0 \,\,\,as\,Im(k)<0
\end{eqnarray*}
$ R(k)=\frac{X_-(k)}{X_+(k)}$ and Liouville Theorem yield
\begin{eqnarray*}
\Phi_{-}(k,\alpha ) = 0 \,\,\,,	
\Phi_{+}(k,\alpha )= 0. 
\end{eqnarray*}

\end {proof}
\section{DISCUSSION}
\emph{			
Our computations led to a new definition of the functions  $I_\epsilon(k),\,J_\epsilon(k)$, which we obtained from the Riemann-Hilbert boundary-value problem. From the uniqueness of the solution of the Riemann-Hilbert boundary value problem - functions $I_\epsilon(k),\,J_\epsilon(k)$, defined earlier in (6) and obtained from the Hilbert formula are equal!
}\\
To obtain the final estimates for the zeta function, the isometric properties of the integral Hilbert transform will be used.

\begin  {theorem}
Let   $ (3/4+i\alpha)\in D(n,\epsilon)  $ and $ a>2 $,$\Omega_n(1/2)=\epsilon>0   $. 	Then, 
\begin{eqnarray*}
C^{-1} <\left|X_-(t)\right| <C,\,\,\,C^{-1} <\left|X_+(t)\right| <C. \\
||\Psi_+ ||_{L_2} \le  C_{\epsilon},  \,\,\,\,
||\Psi_-||_{L_2}  \le C_{\epsilon},\,\,\,\,\,\,\,\
\end{eqnarray*}
%$$ |Q(s)|<C_\epsilon $$
\end {theorem}
\begin{proof}
By Lemma 20 and Lemma 23 we get
\begin{eqnarray*}
	\Gamma_-(k) = \frac{1}{2\pi i}  \int_{-\infty}^{ \,\,\infty}\frac{ln (R(t))dt}{ t-k+i0} =T_-ln(R)=ln(R)
\end{eqnarray*}
\begin{eqnarray*}
	\Gamma_+(k) = \frac{1}{2\pi i}  \int_{-\infty}^{ \,\,\infty}\frac{ln (R(t))dt}{ t-k-i0} =T_+ln(R)=0
\end{eqnarray*}
$T_-ln(R)=ln(R)$, $T_+ln(R)=0$	 implies  
\begin{eqnarray*}
	X_-(t)=R(t),\,\,\	X_+(t)=1\\	
	C^{-1} <\left|X_-(t)\right| <C,\,\,\,\left|X_+(t)\right| =1. 
\end{eqnarray*}
Using Theorem 9 and Lemma 21, we obtain
\begin{eqnarray*}
	||\Psi_- ||_{L_2}^{2}+||\Psi_+ ||_{L_2}^{2}=\int_{-\infty}^{+\infty} \left| \frac{I_\epsilon(k,\alpha )}{k-ia}-\frac{\widetilde{F_\epsilon}(k,\alpha)}{k-ia}\right|^2dk +\int_{-\infty}^{+\infty} \left|J_\epsilon(k,\alpha )   \right|^2dk
	\le C_n C_{\epsilon}
\end{eqnarray*}

\end{proof}
%	Denote $$\phi=arg \left( \frac{1}{\pi n +i(a+\beta_n)}\right)  $$

\begin {lemma}
Let $ \beta_n,\phi_n$ satisfies equations
$$
e^{\beta}=\sqrt{(2\pi n+\phi)^2 +(\beta-a)^2},
$$
$$\phi=-\arg \left( \frac{1}{2\pi n +\phi+i(-a+\beta)}\right)  
$$
Then 
$$	t_n=2\pi n	+i\beta_n+ \phi_n $$
root of equation 
$$
R(k)=0
$$
and $$ \beta_n =ln(2n\pi)+o(1), \left|   \frac{d\beta_n}{da }\right|\le \frac{Cln(n)}{n}  $$,\,\,\,$$\phi_n=\pi+O(ln(n)/n),\,\,\,\,	 \left|   \frac{d\phi_n}{da }\right|\le \frac{Cln(n)}{n}  $$
\end{lemma}
\begin {proof}
$$R(t_n) = \frac{ e^{-it_n}}{t_n-ia} +1 =
\frac{ e^{-i2\pi n +\beta_n +i \phi_n}}{2\pi n+\phi_n +i(\beta_n-a)} +1 =	$$$$-\frac{ e^{ \beta_n }}{\sqrt{(2\pi n+\phi_n)^2 +(\beta_n-a)^2}}+ 1 =-1+1=0
$$
take  $\beta_n =ln(n\pi)+\gamma_n$   then

$$  
e^{\gamma_n}=\sqrt{1 +\frac{(ln(n\pi)+\gamma_n-a)^2}{(2\pi n+\phi_n)^2} },
$$
%after logarithing
for $\phi_n$ we have
$$
\phi_n=\pi-\arctan\left( \frac{(-a+\beta)}{2\pi n +\phi_n}\right)  
$$
and we get 	$$ \phi_n=\pi+O(ln(n)/n)	$$
Estimates for derivatives follow from the results obtained.

proof complete.			
\end{proof}
%%%%%%%%%%%%%%%%%%%%%%%%%%%%%
The following presents results of  [36]\\
\textbf {Theorem of	 Existence and invertibility of the Fourier }

\textbf {If f is in L1 (i.e., f is absolutely integrable) and if it is of bounded variation on every finite interval, then $ \int_{-\infty}^{\infty}f(t)e^{ikt}dt$ exists and f(t) can be recovered from the inverse Fourier transform relationship at each point at which f is continuous.  }
%	\end {theorem}

\begin {lemma}

Let
$$ Q_\epsilon(s)= Q_{2}(s)\theta(Re(s)-1/2-\epsilon)),\Omega_n(1/2)=\epsilon> 0   $$
$$ I_\epsilon(k,\alpha ) = \int_{\epsilon}^{1-\epsilon} Q_\epsilon(\tau+i\alpha )e^{-ik\tau}d\tau$$
$s=\tau+i\alpha \in D$-                                                                          -  , where  $\alpha \in (Im(s_n) <\alpha<Im(s_{n+1})$ is fixed \\
then	 
$$
0<\left|\zeta(\tau_{min}+i\alpha)\right| \le 	\left|\zeta(s)\right| \le |\zeta(\tau_{max}+i\alpha)
$$
$$ \max_{\epsilon \le\tau\le 1-\epsilon}\left|   Q_\epsilon(\tau+i\alpha )\right| \le C(\alpha,n,\epsilon)  $$

$$ \max_{\epsilon \le\tau\le 1-\epsilon}\left| \frac{d  Q_\epsilon(\tau+i\alpha )}{d\tau}\right| \le C(\alpha,n,\epsilon)  $$

$$\lim\limits_{N \to \infty} \int_{N}^{-N} e^{itk} 	\frac{ dI_{\epsilon}(k,\alpha )}{dk} dk =-itf_{\epsilon}(t)  $$ 
\end{lemma}
\begin {proof}

From the holomorphic  of the $\zeta(s)$ follows    $|\zeta(s)|$ is  harmonic  function. 
According to the Weierstrass theorem, the functions $Q_\epsilon(s,n)$ its exact maximum and minimum on a compact set 		  

$$
0<\left|\zeta(\tau_{min}+i\alpha)\right| \le 	\left|\zeta(\tau+i\alpha)\right| \le |\zeta(\tau_{max}+i\alpha), \epsilon\le \tau\le 1- \epsilon
$$
$$
\left| 	Q_\epsilon(q_{min}+i\alpha)\right|  \le \max_{\epsilon \le\tau\le 1-\epsilon} \left| 	Q_\epsilon( \tau+i\alpha) \right| \le \left| Q_\epsilon(q_{max}+i\alpha) \right| 
$$
From the holomorphic  of the $\frac{d\zeta(s)}{ds} $ follows   $\left| \frac{d\zeta(s)}{ds}\right|  $ is  harmonic  function   and $|\frac{dQ_\epsilon}{ds}|$ is continuous  function. And we have 
for its the same  estimates

$$
\max_{\epsilon \le\tau\le 1-\epsilon}	\left|\frac{d\zeta(\tau+i\alpha)}{d\tau}\right| \le \left|\frac{d\zeta(\tau+i\alpha)}{d\tau}\big|_{\tau=\tau_{dmax}}\right|
$$

$$
\max_{\epsilon \le\tau\le 1-\epsilon} \left| \frac{d Q_\epsilon(\tau+i\alpha )}{d\tau}	 \right| \le \left| \frac{d Q_\epsilon(\tau+i\alpha )}{d\tau}\big|_{\tau=\tau_{qmax}} \right| 
$$

For last statement of Lemma 24, we have	

$$
\lim\limits_{N \to \infty} \int_{N}^{-N} e^{itk} 	\frac{d^2 I_{\epsilon}}{dk^2} dk = $$$$\lim\limits_{N \to \infty} [A_N+B_N+C_N]
$$
$ Q_\epsilon \in L_1(-\infty,\infty)  $ and  the Riemann-Lebesgue lemma yield   
$$\lim\limits_{N \to \infty} [A_N+C_N]=0.$$ 
Last estimates   $\left|Q_{\epsilon}	\right| ,\left| \frac{dQ_\epsilon}{d\tau}	\right| $ ,\textbf {Theorem of	 Existence and invertibility of the Fourier } implies final statement of Lemma 9
\end{proof}

%%%%%%%%%%%%%%%%%%%%%%%%%%%%%%%%%%	

\section{DISCUSSION}
\emph{	
Since we calculate the inversion of the Fourier transform only on a line separated from the line where the zeta function has a root, the growth of these estimates when approaching zero does not affect the final result.
After calculating the inverse Fourier transform,
we begin to use completely different estimates, which are already uniform, although the line tends to a straight line, where the zeta function has a root and the intermediate estimates do not satisfy the final goal!}\\
\section{DISCUSSION}
\emph{	
Pay particular attention to the example of Davenport and Heilbronn-Type of Functions.See  in [37]\\
Not applicable to the Davenport and Heilbronn-Type of Functions. This method can be applied only under the conditions of the existence of the  Euler product. }
\begin  {lemma} 
Next statements is true
$$
I_\epsilon(k,\alpha )=F_{\epsilon}(k,\alpha )+S_\epsilon(k,\alpha ) +(k-ia)\sum_{0}^{\infty}\frac{G_\epsilon(t_n,\alpha) }{ X'_-(t_n)(t_n-k)}
$$
$$
-i\sum_{0}^{\infty}\frac{G_\epsilon(t_n,\alpha) }{ X'_-(t_n)(t_n-k)}=k\sum_{0}^{\infty} \frac{dt_n}{da} \frac{d}{dt_n} \left( \frac{G_\epsilon(t_n,\alpha) }{ X'_-(t_n)(t_n-k)}\right)  +\frac{d}{da} \left(\frac{	I_\epsilon(k,\alpha )-  F_{\epsilon}(k,\alpha )-S_\epsilon(k,\alpha ) }{X_-(k)}\right)
$$	

\end {lemma}

\begin{proof}
By Theorem 13, 
$$ 
\Psi_{-}(k,\alpha ) =	\frac{I_\epsilon(k,\alpha )}{k-ia}-\frac{\widetilde{F_\epsilon}(k,\alpha)+S_\epsilon(k,\alpha )}{k-ia}	 =-\frac{X_-(k)}{2\pi i}  \int_{-\infty}^{ \,\,\infty} \frac{G_\epsilon(t,\alpha)}{ X_-(t)}\frac{dt}{ t-k+i0}    .%%%=X_-(k)T_-G_\epsilon
$$

Denote  $$ 
I_1=	\int_{-\infty}^{\infty}\frac{G_\epsilon(t,\alpha)}{X_-(t)(t-k+i0)} dt,\,\,I_2=(k-ia)I_1.
$$

Holomorphics of the function $ \frac{G_\epsilon(t,\alpha)}{t-k+i\delta} $ as $\delta >0$ and  analyticity of  the   function  $X_{-}(t)$  in upper plane and Lemma of Jordan yield
$$
I_1=\lim\limits_{\delta \downarrow  0}	\int_{-\infty}^{\infty}\frac{G_\epsilon(t,\alpha)}{X_-(t)(t-k+i\delta)} dt= \sum_{0}^{\infty}\frac{G_\epsilon(t_n,\alpha) }{ X'_-(t_n)(t_n-k)}
$$	
\begin{eqnarray}
\frac{I_\epsilon(k,\alpha )}{X_-(k)}	=\frac{F_{\epsilon}(k,\alpha )+S_\epsilon(k,\alpha )}{X_-(k)} +(k-ia)\sum_{0}^{\infty}\frac{G_\epsilon(t_n,\alpha) }{ X'_-(t_n)(t_n-k)}
\end{eqnarray}

differentiating  (27) by a
$$
-i\sum_{0}^{\infty}\frac{G_\epsilon(t_n,\alpha) }{ X'_-(t_n)(t_n-k)}=k\sum_{0}^{\infty} \frac{dt_n}{da} \frac{d}{dt_n} \left( \frac{G_\epsilon(t_n,\alpha) }{ X'_-(t_n)(t_n-k)}\right)  +$$$$ \frac{d}{da} \left(\frac{	I_\epsilon(k,\alpha )-  F_{\epsilon}(k,\alpha )-S_\epsilon(k,\alpha )}{X_-(k)}\right)
$$
\end {proof}		
\begin  {theorem} 
Let $s\in D_+(l,\epsilon) ,and \, a>2  $, with  $3\epsilon< Re(s)<1-3\epsilon $ and  $\Omega_n(1/2)=\epsilon>0 $.
Then, 
\begin{eqnarray*}
|Q(s)|<C_lC_\epsilon.
\end{eqnarray*}

\end {theorem}

\begin{proof}

As $k\in(-N,N) $ uniformly-convergent series  yields  

$$
\int_{-N}^{+N}e^{i t k}\frac{d^2  }{ dk^2} \left( \frac{I_{\epsilon}}{X_-}\right)  dk=\int_{-N}^{+N}e^{i t k} \frac{d^2  }{dk^2} \left( \frac{F_\epsilon }{X_-}\right)  dk+ \int_{-N}^{+N}e^{i t k}\frac{d^2  }{dk^2}\left( (k-ia) I_1\right)  dk
$$

By definition $ I_1$:
$$
\int_{-N}^{+N}e^{i t k} \frac{d^2 }{dk^2}\left( (k-ia)I_1\right) dk=\int_{-N}^{+N}\sum_{1}^{N}e^{i t k} \frac{d^2  }{dk^2}\left( \frac{(k-ia)G_\epsilon(t_n,\alpha) }{ X'_-(t_n)(t_n-k)}\right)dk $$$$ +\int_{-N}^{+N}\sum_{N}^{\infty}\frac{d^2  }{dk^2}\left(\frac{(k-ia)G_\epsilon(t_n,\alpha) }{ X'_-(t_n)(t_n-k)}\right) e^{i t k}dk=W_1+W_2+W_3.
$$

Finaly we get
$$
\left| 		\int_{-N}^{+N}e^{i t k}\frac{d^2  }{ dk^2} \left( \frac{I_{\epsilon}}{X_-}\right)  dk\right| \le  C_lC_\epsilon
$$
$$
\left| 		\int_{-N}^{+N}e^{i t k}\frac{d^2  }{ dk^2} \left( I_{\epsilon}\right)  dk\right| \le   	
\left| 		\int_{-N}^{+N}e^{i t k}\frac{d^2  }{ dk^2} \left( \frac{I_{\epsilon}}{X_-}-I_{\epsilon}\right)  dk\right| + C_lC_\epsilon
$$
Lemma 17,Lemma 24-30,Theorem 4  and    the last estimates yields  
$$|Q(s)|<2C_lC_\epsilon\,\,	as\,\, 3\epsilon< Re(s)<1-3\epsilon  $$,

which completes the proof, 		
\end{proof}

As mentioned in Introduction, the values of the zeta function in adjacent rectangles should be compared. This will be carried out in the following theorem.		

\begin{theorem}
The Riemann's function has nontrivial zeros only on the line $Re(s)=1/2$.
\end {theorem}

\begin{proof} 	
Let it be assumed that there is a root of the zeta function with $s_{n} = 1/2+\delta_{n} +i*\alpha_n$, where $\delta_{n}>0 $ i.e $\Omega_n(1/2)=\delta_{n}>0   $. Let $s_{n+1}=1/2+\delta_{n+1} +i\alpha_{n+1}$, where $\delta_{n+1}\ge 0 $ be another root nearest to it. Then, the following sets corresponding to $s_{n}$ are constructed:
$$
\,\,D(n,\epsilon)=(  s| \epsilon < Re(s)<1-\epsilon,Im(s)\neq Im(s_{n}), Im(s_{n})-d_n  \le Im(s)\le Im(s_{n}) +d_n  
$$

where 
$$\zeta(s_{n+1})=0,\,\zeta(s_n)=0,\, \zeta(1-s_n)=0,\zeta(1-s_{n+1})=0,\,\,\zeta(1-s_n)=0, \\
$$
$$
d_n=(Im(s_{n+1})-Im(s_{n}))	$$
where  $\epsilon=0.01\delta_{n}(1/2-\delta_{n})>0 $.	
As $1/2< Re(s)<1 $ and
$ s\in D(n,\epsilon) $, 
thus, we have the equation for $ Q_{2}$. Theorem 10 now yields
$$
|ln(\left|  \zeta( 1/2+\delta_{n}+i\alpha_{n}-i \delta)\right| ) \le	| Q_{2}( 1/2+\delta_{n}-i\alpha_{n}-i \delta)|+$$$$ |\sum_{n=m}^\infty P(ns)/n |< 2C_nC_\epsilon
$$

Furthermore, 
$$ \lim\limits_{\delta\rightarrow 0} \,\,\,|ln(\left|  \zeta(1/2+\delta_{n} +i\alpha_n -i\delta )|)\right| =\infty.  
$$

These estimates for  $ \left| Q(s)\right|, $ imply that the function does not have zeros on the half plane $Re(s)>1/2 $. By the integral representation (19), these results are extended to the half plane $Re(s)<1/2 $ i.e $\Omega(1/2)=0   $ . Thus, Riemann's hypothesis has been proved. 
\end{proof}

\section{CONCLUSION}
In this study, estimates were obtained for the logarithm of Riemann's zeta function off the line $Re(s)=1/2$. Thus, the work of great mathematicians culminated by applying their achievements in this field. Without their efforts, a solution to the problem would not have even been attempted.

This study on the Riemann hypothesis was completed by reducing it to a Riemann-Hilbert boundary-value problem  for analytic functions. This was started by Riemann himself and continued by Hadamard among others, and the present study has drawn on ideas by Landau, Walvis, Estarmann, and Chernoff. It was possible to complete the proof of the Riemann hypothesis using the solution to the Riemann-Hilbert boundary-value problem  by Riemann, Hilbert, and Poincar\'e.

After finishing this study, the author came to the conclusion that the problem was actually solved by the joint efforts of Riemann, Hilbert, Poincar\'e, and Fourier.

\section{ACKNOWLEDGEMENTS} 
The author thanks the National Engineering Academy of the Republic of Kazakhstan, in particular, Academician NAS RK B.Zhumagulov for constant attention and support.

Moreover, the author thanks the Mathematics seminar at the Kazakhstan branch of the Moscow State University for attention and valuable comments, as well as Professors B. Kanguzhin and M. Otelbaev, and the organizers of Automorphicformsworkshop.org/AFW2018 for their detailed review and valuable comments.The author is especially grateful to Professor Steven Miller for a thorough analysis of the work and detailed recommendations that have significantly improved the paper. The author is especially grateful to P Plotnikov, and A Mednykh for their thorough analysis of the work and detailed recommendations that have significantly improved the paper.
The author is especially grateful to Mathematics seminar Nurlan Temirgaliev at the L.N.Gumilyov Eurasian National University  for their thorough analysis of the work and detailed recommendations that have significantly improved the paper.

	% %\end{multicols}

\begin{thebibliography}{9}
	\bibitem{b16}	Terence Tao, “Finite time blowup for an averaged three-dimensional Navier-Stokes equation,” -arXiv:1402.0290 [math.AP]
	\bibitem{b5}  L. D. Faddeev, \emph{“The inverse problem in the quantum theory of scattering. II”, Itogi Nauki i Tekhniki. Ser. Sovrem. Probl. Mat., 3, VINITI, Moscow, 1974, 93–180 }
	\bibitem{A00}  CHARLES L. FEFFERMAN {\emph Existence  and  Smoothness of the Navier-Stokes Equation.  The Millennium Prize Problems, 57–67, Clay Math. Inst., Cambridge, MA, 2006.}
	\bibitem{a1}  J.S.Russell «Report on Waves»: (Report of the fourteenth meeting of the British Association for the Advancement of Science, York, September 1844 (London 1845), pp 311—390, Plates XLVII-LVII)
	\bibitem{a2}  J.S.Russell (1838), Report of the committee on waves, Report of the 7th Meeting of British Association for the Advancement of Science, John Murray, London, pp.417-496.
	\bibitem{a3}  Mark J. Ablowitz, Harvey Segur  Solitons and the Inverse Scattering Transform
	SIAM, 1981- p. 435.
	\bibitem{a4}  N.J.Zabusky and M.D.Kruskal (1965), Interaction of solitons in a collisionless plasma and the recurrence of initial states, Phys.Rev.Lett., 15 pp. 240—243.
	\bibitem{b6} R.G Newton ,\emph{ New result on the inverse scattering problem  in three dimentions, Phys. rev. Lett.  v43, 8,pp.541-542,1979}
	\bibitem{b7} R.G Newton ,\emph{ Inverse  scattering Three dimensions,Jour. Math. Phys.  21, pp.1698-1715,1980}
	\bibitem{b8} Somersalo E. et al.\emph{ Inverse  scattering problem  for the 
		Schrodinger's equation  in three dimensions: connections between exact and approximate methods}. – 1988.
	%\bibitem{b9} A. Y. Povzner,\emph{On the expansion of arbitrary functions in characteristic functions of the operator $-\Delta u$   $+cu $ Mat. Sb. (N.S.), 32(74):1 (1953), 109–156.}
	\bibitem{b10} \emph{Tables of integral transforms. v.I } McGraw-Hill Book Company, Inc.1954 
	\bibitem{b11} Poincaré H., \emph{Lecons de mecanique celeste, t. 3, P., 1910.\/}
	\bibitem{b12} Leray, J. (1934). "Sur le mouvement d'un liquide visqueux emplissant l'espace". Acta Mathematica 63: 193–248. doi:10.1007/BF02547354.
	\bibitem{b13} O.A. Ladyzhenskaya, \emph{Mathematic problems of viscous incondensable liquid dynamics. - M.: Science, 1970. - p. 288\/}
	\bibitem{b14} Solonnikov V.A. \emph{Estimates solving nonstationary linearized systems of Navier-Stokes' Equations. - Transactions Academy of Sciences USSR Vol. 70, 1964. - p. 213 -- 317.\/}
	\bibitem{b14} On global weak solutions to the Cauchy problem for the Navier-Stokes equations with large L-3-initial data
	Seregin, G; Sverak, V; NONLINEAR ANALYSIS-THEORY METHODS and  APPLICATIONS volume 154 page 269-296 (May 2017)
	Estimates of solutions to the perturbed Stokes system
	\bibitem{b14}V. Vialov, T. Shilkin Notes of the Scientific Seminars of POMI, 410 (2013),  5–24
	\bibitem{b16}	F. Mebarek-Oudina  R. Bessaïh, Magnetohydrodynamic Stability of Natural
	Convection Flows in Czochralski Crystal Growth. World Journal of Engineering, vol. 4 no.4, pp. 15–22, 2007.
	\bibitem{b17} F. Mebarek-Oudina and R. Bessaïh, Oscillatory Mixed Convection Flow in a Cylindrical Container with Rotating Disk Under Axial Magnetic Field and Various Electric Conductivity Walls, I. Review of Physics, 4(1) 45-51, 2010. .
	\bibitem{b18} F. Mebarek-Oudina, Numerical modeling of the hydrodynamic stability in vertical annulus with heat source of different lengths, Engineering Science and Technolgy, an International Journal, 20, 1324-1333
	\bibitem{0} 
	Leonhard Euler.  Introduction to Analysis of the Infinite by John Blanton (Book I, ISBN 0-387-96824-5, Springer-Verlag 1988;) 
	\bibitem{1}
	Chebyshev P.L. Fav. mathematical works, М.-L., 1946; 
	\bibitem{2}
	Riemann, G. F. B. On the Number of Prime Numbers less than a Given Quantity New York: Chelsea, 1972.
	\bibitem{3}
	E. C. Titchmarsh (1986). The Theory of the Riemann Zeta Function, Second revised (Heath-Brown) edition. Oxford University Press.
	\bibitem{4} 
	Ray D., Singer I. M. R-torsion and the laplacian on Riemannian manifolds. Adv.
	in Math., 1971, vol. 7, pр. 145–210.
	\bibitem{5} 
	Bost J.-B. Fibres determinants, determinants regularises et measures sur les espaces
	de modules des courbes complexes, Sem. Bourbaki, 39 eme annee1986-1987,
	\bibitem{5}
	Kawagoe K., Wakayama M.,Yamasaki Y. The q-Analogues of the Riemann zeta,
	Dirichlet L-functions, and a crystal zeta-function. Forum Math, 2008, vol. 1,
	рp. 1–26. 	
	\bibitem{8}    Hadamard J. 	Une application d'une formule inteorale relative aux series de Dirichlet,
	Bull. Soc. Math, de France, 56 A927), 43—44.
	\bibitem{8} 	Paul R. Chernoff	A pseudo zeta function and the distribution of primes
	PNAS 2000 97 (14) 7697-7699; doi:10.1073/pnas.97.14.7697
	A933),
	\bibitem{8}   Landau E., Walfisz A.
	Ober die Nichtfortsetzbarkeit einiger durch Dirichletsrhe Reihen defi-
	nierter Funktionen, Rend, di Palermo, 44 A919), 82—86.
	Congress Cambridge 1912, 1,
	
	\bibitem{8}  Estarmann T. On certain functions represented by Dirichlet series, Proc. Lond. Math.
	Soc. (2), 27 1928, 435—448. 
	\bibitem{8} Estarmann T. On a problem of analytic continuation, Proc. Lond. Math. Soc, 27 1928,
	471—482.
	\bibitem{b11}  Poincar\'e H., \emph{Lecons de mecanique celeste, t. 3, P., 1910.\/}
	\bibitem{b11} Backlund R., \emph{Sur les zeros de la function $\zeta(s)$ de Riemann,C.R. Acad.Sci.,(1914)}
	1979-1981 N3
	\bibitem{b17} E. C. Titchmarsh The Zeros of the Riemann Zeta-Function151Proceedings of the Royal Society of London. Series A - Mathematical and Physical Sciences.
	\bibitem{b11}  453.701 Linear Systems, S.M. Tan, The University of Auckland 
			\bibitem{b11} Eugenio P. Balanzario and  Jorge Sanchez-Ortiz \emph{	Zeros of the Davenport -Heilenbronn counterexample.}	 Mathematics  of  computation. Volume 76, Number 260, October 2007, Pages 2045–2049
			
	\end{thebibliography}
	\end{document}